\documentclass[11pt,a4paper,oneside,reqno]{amsart}
\usepackage{amsfonts}
\usepackage{amssymb}
\usepackage{mathrsfs}
\usepackage[super]{nth}
\usepackage{amsthm}
\usepackage{soul}
\usepackage[defaultcolor=orange]{changes}
\usepackage{amsmath}
\usepackage{enumitem}
\usepackage{mathtools}
\usepackage{adjustbox}
\usepackage[bottom]{footmisc}
\usepackage{pifont}
\usepackage{float}
\usepackage{pgfplots}
\usepackage{color}
\usepackage{url}
\usepackage{tikz-cd}
\usepackage{graphicx}
\usepackage{csquotes}
\usepackage[top=1in, bottom=1in]{geometry}
\graphicspath{ {images/} }

\usepackage[hidelinks]{hyperref}
\usepackage{makecell}
\usepackage{fancyhdr}

\newcommand\cyr
{
\renewcommand\rmdefault{wncyr}
\renewcommand\sfdefault{wncyss}
\renewcommand\encodingdefault{OT2}
\normalfont
\selectfont
}
\DeclareTextFontCommand{\textcyr}{\cyr}

\newlist{case}{enumerate}{1}
\setlist[case,1]{
  label={\textsc{Case}~\arabic*:},
  leftmargin=*,
  align=left,
  labelsep=5mm
}

\newlist{step}{enumerate}{1}
\setlist[step,1]{
  label={\textsc{Step}~\arabic*:},
  leftmargin=*,
  align=left,
  labelsep=5mm
}

\allowdisplaybreaks
\usepackage[OT2,OT1]{fontenc}
\usepackage{wasysym}
\DeclareSymbolFont{extraup}{U}{zavm}{m}{n}
\DeclareMathSymbol{\varheart}{\mathalpha}{extraup}{86}
\DeclareMathSymbol{\vardiamond}{\mathalpha}{extraup}{87}
\DeclareMathSymbol{\varspade}{\mathalpha}{extraup}{85}
\pgfdeclarelayer{nodelayer}
\pgfdeclarelayer{edgelayer}
\pgfsetlayers{nodelayer,edgelayer,main}
\def \proof {\noindent\textit{Proof.}\hspace{5mm}}
\def \proofsketch {\noindent\textit{Proof (Sketch).}\hspace{5mm}}
\def \bs {$\hfill \blacksquare$\\}
\def \bsn {$\hfill \blacksquare$}
\def \sq {$\hfill \square$\\}
\def \sqn {$\hfill \square$}

\DeclareMathOperator{\Inv}{Inv}
\DeclareMathOperator{\Aut}{Aut}

\DeclareMathOperator{\Hom}{Hom}

\DeclareMathOperator{\Subg}{Subg}

\DeclareMathOperator{\Gal}{Gal}

\DeclareMathOperator{\Z}{\mathbb{Z}}
\DeclareMathOperator{\N}{\mathbb{N}}
\DeclareMathOperator{\Q}{\mathbb{Q}}

\DeclareMathOperator{\R}{\mathbb{R}}
\DeclareMathOperator{\C}{\mathbb{C}}

\def\forkindep{\mathrel{\raise0.2ex\hbox{\ooalign{\hidewidth$\vert$\hidewidth\cr\raise-0.9ex\hbox{$\smile$}}}}}
\definecolor{supergreen}{RGB}{0, 170, 0}
\definecolor{superred}{RGB}{170, 0, 0}
\definecolor{darkred}{RGB}{100, 0, 0}

\newcommand{\Th}{\mbox{\ttfamily Th}}

\newcommand{\TTh}{\mbox{\emph{\ttfamily Th}}}

\newcommand{\Kbad}{K_{\mbox{\tiny bad}}}

\newcommand{\iso}{\cong}
\newcommand{\defeq}{\mathrel{\mathop:}=}

\newtheorem{thm}{Theorem}[section]
\newtheorem{prop}[thm]{Proposition}
\newtheorem{cor}[thm]{Corollary}
\newtheorem{lem}[thm]{Lemma}

\newtheorem{prob}[thm]{Problem}

\theoremstyle{definition}
\newtheorem{definition}[thm]{Definition}
\newtheorem{cons}[thm]{Construction}
\newtheorem{exmp}[thm]{Example}
\newtheorem{remark}[thm]{Remark}
\newtheorem{innercustomthm}{Theorem}
\newenvironment{customthm}[1]
  {\renewcommand\theinnercustomthm{#1}\innercustomthm}
  {\endinnercustomthm}

\newtheorem{innercustomconj}{Conjecture}
\newenvironment{customconj}[1]
  {\renewcommand\theinnercustomconj{#1}\innercustomconj}
  {\endinnercustomconj}

\pgfplotsset{compat=1.14}

\begin{document}
\title{Finite Undecidability in Fields II: PAC, PRC \& P$p$C Fields}
\author{Brian Tyrrell}
\address{Mathematical Institute, Woodstock Road, Oxford OX2 6GG.}
\email{brian.tyrrell@maths.ox.ac.uk}

\begin{abstract}
A field $K$ in a ring language $\mathcal{L}$ is \emph{finitely undecidable} if $\mbox{Cons}(\Sigma)$ is undecidable for every nonempty finite $\Sigma \subseteq \mbox{\ttfamily Th}(K; \mathcal{L})$. We adapt arguments originating with Cherlin-van den Dries-Macintyre/Ershov (for PAC fields) and Haran (for PRC fields) to prove all PAC and PRC fields are finitely undecidable. We describe the difficulties that arise in adapting the proof to P$p$C fields, and show no (\emph{bounded}) P$p$C field is finitely axiomatisable. This work is drawn from the author's PhD thesis \cite[Chapter 4]{btyrrelthesis} and is a sequel to \cite{btyrrel}.
\end{abstract}

\maketitle
\vspace{-3mm}
\section{Introduction}
\setul{1.5pt}{.4pt}

\noindent The author was motivated to consider this topic by the following question:

\begin{prob}\label{1}
Does there exist an infinite, finitely axiomatisable field? 
\end{prob}

This (open) problem was posed explicitly by I.\ Kaplan at the 2016 Oberwolfach workshop on \emph{Definability and Decidability Problems in Number Theory} \cite[Q4]{ober}, though existed as folklore before this. One approach to \emph{Problem \ref{1}} was established by Ziegler in 1982 \cite{ziegler}, and generalised further by Shlapentokh \&\ Videla \cite{shlapentokhvidela}. We forward the following definition\footnote{Shlapentokh \& Videla \cite{shlapentokhvidela} call this property \emph{finite hereditary undecidability}; for notational ease we remove the word ``hereditary''.}:

\begin{definition}
A theory $T$ in a language $\mathcal{L}$ is \emph{finitely undecidable} if every finitely axiomatised $\mathcal{L}$-subtheory of $T$ is undecidable. (An $\mathcal{L}$-structure is \textit{finitely undecidable} if its $\mathcal{L}$-theory is.)
\end{definition}

Ziegler in 1982 proved $\C$, $\widetilde{\mathbb{F}_p(t)}$, $\R$, and $\Q_p$ are finitely undecidable (in the language of rings) \cite[Folgerung, p.\ 270]{ziegler}, and the author recently extended these results to henselian NIP nontrivially valued fields \cite[Corollary 5.12]{btyrrel} (for undecidability in the language of valued fields). For this paper, we are motivated to consider what other finitely undecidable fields exist -- we will show that \emph{PAC and PRC fields are finitely undecidable} in the language of rings. Model-theoretically, this is in a strictly different direction to \cite{btyrrel}, as there we focused on NIP theories, and here (we shall see) the focus is on \emph{simple} theories (and generalisations of simplicity). Recall the definition of a \emph{simple} theory by Shelah \cite{shelahntp}, expounded upon in \cite{wagnersimpletheories}.

In \S \ref{sect31} we follow a suggestion of E.\ Hrushovski to adapt the arguments of Cherlin, van den Dries \&\ Macintyre \cite{cvddmac} and independently Ershov \cite{ershov} to prove \emph{every PAC field is finitely undecidable} (\emph{Corollary \ref{mainpacfieldundeccor}}). In \S \ref{theprcfields} we adapt the work of Haran \& Jarden \cite{haran, haranjarden} to show more generally \emph{every PRC field is finitely undecidable} too (\emph{Corollary \ref{mainprcfieldundeccor}}). As a consequence we answer two open questions of Shlapentokh \& Videla \cite[\S 6]{shlapentokhvidela} -- see \emph{Remark \ref{theopenqs}}. We are unable to use this method to prove finite undecidability of P$p$C fields -- \S \ref{stuffforppc}
explains why -- but we can determine no \textit{bounded} P$p$C field is finitely axiomatizable (\textit{Theorem \ref{noppcfinax}}) after setting up the correct machinery to adjust work of Haran \& Jarden \cite{haranjardenppc}. \underline{Familiarity with \cite{cvddmac}} (for which Chatzidakis gives an excellent overview, in \cite[Appendix 1]{chatzidakis}) \underline{and \cite{haranjarden, haranjardenppc} is assumed throughout}.

This gives the author confidence to forward the following problem:

\begin{prob}\label{3}
Does there exist an infinite field that is \emph{not} finitely undecidable?
\end{prob}

If \emph{Problem \ref{3}} is resolved in the negative, this of course gives a negative answer to \emph{Problem \ref{1}}. Assuming powerful classification-theoretic conjectures, progress can be made on at least understanding the shape this problem takes: model-theoretically `tame' structures will be finitely undecidable, and already many model-theoretically `wild' structures have this property (due to their interpretation of arithmetic; e.g.\ the theories of all number fields and global function fields by J.\ Robinson \cite{rob59}, resp.\ Rumely \cite{rumely}, and more generally the theory of any positive characteristic function field by Eisentr\"{a}ger \&\ Shlapentokh \cite{eisenshlap} or any infinite finitely generated field by Poonen \cite[Remark 5.2]{poonenpaper}).

\section{Preliminaries on Undecidability}

All of our undecidability results rely on a theorem of Ershov, which we give below. See \cite[Chapter 5, \S 1]{ershovrussian} in comparison. First, some definitions:

\begin{definition}\label{interp}
\cite[p.\ 212 \& \S5.3 Remark 4]{hodges}. Let $\mathcal{L}_0$, $\mathcal{L}_1$ be languages, $A$ an $\mathcal{L}_1$-structure and $B$ an $\mathcal{L}_0$-structure. We define an \textit{interpretation of $B$ in $A$} to be:
\begin{enumerate}
    \item An $\mathcal{L}_1$-formula $\delta(x_1, \dots, x_n)$;
    \item For every atomic $\mathcal{L}_0$-formula $\phi(y_1, \dots, y_m)$, an $\mathcal{L}_1$-formula $\phi'(\overline{x}_1, \dots, \overline{x}_m)$ in which the $\overline{x}_i$ are disjoint $n$-tuples of distinct variables;
    \item A surjective function $f : \delta(A^n) \rightarrow B$,
\end{enumerate}
such that for all atomic $\mathcal{L}_0$-formulae $\phi$ and $\overline{a}_i \in \delta(A^n)$,
\begin{equation}
B \models \phi(f(\overline{a}_1), \dots, f(\overline{a}_m)) \iff A \models \phi'(\overline{a}_1, \dots, \overline{a}_m).\tag{$\dagger$}
\end{equation}
We say \textit{$B$ is interpretable in $A$} if there exists a interpretation of $B$ in $A$. We say $B$ is interpretable in $A$ \emph{with parameters} if there exists a set $S \subset A$ such that $B$ is interpretable in the $\mathcal{L}_1(S)$-structure $A$.
\end{definition}

\begin{definition}\label{uniformparam}
Let $\mathcal{L}_0, \mathcal{L}_1$ be finite languages, $K_0$ a class of $\mathcal{L}_0$-structures and $K_1$ a class of $\mathcal{L}_1$-structures. We say \emph{$K_0$ is uniformly interpretable in $K_1$} if there exists:
\begin{enumerate}
    \item An $\mathcal{L}_1$-formula $\delta(x_1, \dots, x_n)$;
    \item For every atomic $\mathcal{L}_0$-formula $\phi(y_1, \dots, y_m)$, an $\mathcal{L}_1$-formula $\phi'(\overline{x}_1, \dots, \overline{x}_m)$ in which the $\overline{x}_i$ are disjoint $n$-tuples of distinct variables;
\end{enumerate}
such that for any $B \in K_0$ there exists:
\begin{itemize}
    \item $A \in K_1$;
    \item A surjective function $f_{AB} : \delta(A^n) \rightarrow B$ such that for all atomic $\mathcal{L}_0$-formulae $\phi$ and $\overline{a}_i \in \delta(A^n)$, $B \models \phi(f_{AB}(\overline{a}_1), \dots, f_{AB}(\overline{a}_m)) \Leftrightarrow A \models \phi'(\overline{a}_1, \dots, \overline{a}_m)$.
\end{itemize}
We say $K_0$ is \emph{uniformly interpretable with parameters} in $K_1$ if there exists a finite expansion $\mathcal{L}_1'$ of $\mathcal{L}_1$ by constant symbols such that $K_0$ is uniformly interpretable in $K_1$ as a class of $\mathcal{L}_1'$-structures.
\end{definition}

\begin{remark}\label{rema}
A class $K_0$ of $\mathcal{L}_0$-structures being uniformly interpretable with parameters in a class $K_1$ of $\mathcal{L}_1$-structures is what Ershov refers to as `$K_0$ being \emph{relatively elementarily definable} in $K_1$', in \cite[pp.\ 271--272]{ershovrussian}. Ershov notes if $K_1^*$ is a class of $\mathcal{L}_1$-structures containing the class $K_1$, and $K_0$ is uniformly interpretable with parameters in $K_1$, then $K_0$ is uniformly interpretable with parameters in $K_1^*$ (\cite[p.\ 272]{ershovrussian}). \sq
\end{remark}

For the remainder of this section, all languages will be finite, hence all interpretations \textit{recursive}. The following definitions and their presentation were very gratefully suggested to the author by E.\ Hrushovski (see \cite[pp.\ 221--222]{hodges} for discussion):

\begin{definition}
Let $\mathcal{L}_0, \mathcal{L}_1$ be finite languages, $K_0$ a class of $\mathcal{L}_0$-structures and $K_1$ a class of $\mathcal{L}_1$-structures. We say $K_0$ is uniformly interpretable \emph{in the strict sense} in $K_1$ if $K_0$ is uniformly interpretable in $K_1$, and (following the notation of \emph{Definition \ref{uniformparam}}) for every $A \in K_1$ there exists $B \in K_0$ and a surjective function $f_{AB} : \delta(A^n) \rightarrow B$ such that for all atomic $\mathcal{L}_0$-formulae $\phi$ and $\overline{a}_i \in \delta(A^n)$, 
$$B \models \phi(f_{AB}(\overline{a}_1), \dots, f_{AB}(\overline{a}_m)) \iff A \models \phi'(\overline{a}_1, \dots, \overline{a}_m).$$
\end{definition}

Let $T$ be a theory in a language $\mathcal{L}$. Denote by $\mathbb{K}_T$ the class of all $\mathcal{L}$-structures satisfying $T$, and if $\mathbb{K}$ is a class of $\mathcal{L}$-structures, denote by $\Th(\mathbb{K}; \mathcal{L})$ the common $\mathcal{L}$-theory of $M \in \mathbb{K}$.

\begin{definition}\label{a11}
Let $\mathcal{L}_0, \mathcal{L}_1$ be finite languages, $T_0$ an $\mathcal{L}_0$-theory and $T_1$ an $\mathcal{L}_1$-theory. We say \emph{$T_0$ is interpretable} (resp.\ \emph{interpretable with parameters}, resp.\ \emph{interpretable in the strict sense}) \emph{in $T_1$} if $\mathbb{K}_{T_0}$ is uniformly interpretable (resp.\ uniformly interpretable with parameters, resp.\ uniformly interpretable in the strict sense) in $\mathbb{K}_{T_1}$. 
\end{definition}

\noindent This definition leads to a nice transfer of undecidability from $T_0$ to $T_1$:

\begin{lem}\label{a12}
Let $\mathcal{L}_0, \mathcal{L}_1$ be finite languages, $T_0$ an $\mathcal{L}_0$-theory and $T_1$ an $\mathcal{L}_1$-theory. Suppose $T_0$ is interpretable in the strict sense in $T_1$. If $T_0$ is undecidable, then $T_1$ is undecidable.
\end{lem}

\proof
We will argue that if $T_1$ is decidable, so too is $T_0$. For $\varphi \in \mbox{Sent}(\mathcal{L}_0)$, $\varphi \in T_0$ if and only if for all $M \in \mathbb{K}_{T_0}$, $M \models \varphi$. By assumption, for $M \in \mathbb{K}_{T_0}$ there exists $N \in \mathbb{K}_{T_1}$ such that $N \models \varphi'$. Furthermore, for all $N \in \mathbb{K}_{T_1}$, $N \models \varphi'$: otherwise $N \models \neg \varphi'$, and as $\neg\varphi' = (\neg \varphi)'$, by \emph{Definition \ref{a11}} there exists $M \in \mathbb{K}_{T_0}$ with $M \models \neg\varphi$, a contradiction. We conclude $\varphi \in T_0 \iff \varphi' \in T_1$. As the reduction map $-'$ is recursive and $T_1$ is decidable, so too is $T_0$. \bs

\begin{thm}\label{a4}
\emph{\textbf{(Ershov).}} Let $\mathcal{L}_0, \mathcal{L}_1$ be finite languages, $T_0$ an $\mathcal{L}_0$-theory and $T_1$ an $\mathcal{L}_1$-theory. Suppose $T_0$ is interpretable with parameters in $T_1$; if $T_0$ is hereditarily undecidable, then $T_1$ is hereditarily undecidable.
\end{thm}

This result is originally due to Ershov (\cite[Chapter 5, \S 1.4, Theorem 2]{ershovrussian}; cf.\ \emph{Remark \ref{rema}} for the terminology ``relatively elementarily definable'').\\

\proof
By \emph{Definitions \ref{uniformparam} \& \ref{a11}}, for some finite expansion by constants $\mathcal{L}_1^{\flat} \supseteq \mathcal{L}_1$, $T_0$ is interpretable in $T_1$ as an $\mathcal{L}_1^{\flat}$-theory. Consider the reduction map $-' : \mbox{Form}(\mathcal{L}_0) \rightarrow \mbox{Form}(\mathcal{L}^{\flat}_1)$ from \emph{Definition \ref{uniformparam}}. We claim $T^{\circ} = \{\varphi \in \mbox{Sent}(\mathcal{L}_0) \mbox{ : } \varphi' \in T_1\}$ is a subtheory of $T_0$. Indeed, define $\mathbb{K}_{T_0}^{\circ}$ to be the class of $\mathcal{L}_0$-structures $M^{\circ}$ such that there exists $N \in \mathbb{K}_{T_1}$ and a surjective function $f_{NM^{\circ}} : \delta(N^n) \rightarrow M^{\circ}$ such that for all atomic $\mathcal{L}_0$-formulae $\phi(x_1, \dots, x_m)$ and $\overline{a}_i \in \delta(N^n)$,
$$M^{\circ} \models \phi(f_{NM^{\circ}}(\overline{a}_1), \dots, f_{NM^{\circ}}(\overline{a}_m)) \iff N \models \phi'(\overline{a}_1, \dots, \overline{a}_m).$$

By assumption, $\mathbb{K}_{T_0} \subseteq \mathbb{K}_{T_0}^{\circ}$, hence $\Th(\mathbb{K}_{T_0}^{\circ}; \mathcal{L}_0) \subseteq T_0$, and $\mathbb{K}_{T_0}^{\circ}$ is uniformly interpretable in the strict sense in (the class of $\mathcal{L}_1^{\flat}$-structures) $\mathbb{K}_{T_1}$. By definition, for $\varphi \in \mbox{Sent}(\mathcal{L}_0)$, $\varphi \in \Th(\mathbb{K}_{T_0}^{\circ}; \mathcal{L}_0) \iff \varphi' \in T_1$ and $T^{\circ} = \Th(\mathbb{K}_{T_0}^{\circ}; \mathcal{L}_0)$ is uniformly interpretable in the strict sense in $T_1$. We conclude the undecidability of $T_1$ from \emph{Lemma \ref{a12}}. 

If $S$ is a subtheory of $T_1$, $\mathbb{K}_{T_1} \subseteq \mathbb{K}_S$. By \emph{Remark \ref{rema}}, $T_0$ is interpretable with parameters in $S$, hence the above proof applies. We conclude $T_1$ is \textit{hereditarily} undecidable.

Finally, we claim that $T_1$ is hereditarily undecidable as an $\mathcal{L}_1$-theory. (We will follow \cite[Proposition 11.2]{shoenfield} for this.) First note that if $\phi \in \mbox{Sent}(\mathcal{L}_1^{\flat})$, there exists $\phi_{\flat}(x_1,\dots, x_k) \in \mbox{Form}(\mathcal{L}_1)$ such that $\phi$ is obtained from $\phi_{\flat}$ by replacing the free occurrences of $x_1, \dots, x_k$ by constants $c_1, \dots, c_k \in \mathcal{L}_1^{\flat} \setminus \mathcal{L}_1$. Now, let $S$ be an $\mathcal{L}_1$-subtheory of $T_1$; i.e.\ $S$ is a subtheory of $\Th(\mathbb{K}_{T_1}|_{\mathcal{L}_1}; \mathcal{L}_1)$, where $\mathbb{K}_{T_1}|_{\mathcal{L}_1}$ is the class of $\mathcal{L}_1$-structures $M$ such that $M \in \mathbb{K}_{T_1}|_{\mathcal{L}_1} \iff M^+ \in \mathbb{K}_{T_1}$, where $M^+$ is an expansion of $M$ to $\mathcal{L}_1^{\flat}$. Let $S^{\flat} = \mbox{Cons}(S)$ as a subtheory of the $\mathcal{L}_1^{\flat}$-theory $T_1$. We claim $\phi \in S^{\flat} \iff \forall x_1, \dots, x_k\, \phi_{\flat}(x_1, \dots, x_k) \in S$.

Indeed, $\forall x_1, \dots, x_k\, \phi_{\flat}(x_1, \dots, x_k) \in S \implies \phi \in S^{\flat}$ is immediate. Suppose $\phi \in S^{\flat}$, and consider $N \in \mathbb{K}_S$ with $n_1, \dots, n_k \in N$. Expand $N$ to an $\mathcal{L}_1^{\flat}$-structure $N^+$ by setting $c_1^{N^+}\!\! = n_1, \dots, c_k^{N^+}\!\! = n_k$. As $N \models S$, $N^+ \models S^{\flat}$, hence $N^+ \models \phi$ and thus $N \models \phi_{\flat}(n_1, \dots, n_k)$. As $N \in \mathbb{K}_S$, $n_1, \dots, n_k \in N$ were arbitrary, $\forall x_1, \dots, x_k\, \phi_{\flat}(x_1, \dots, x_k) \in \Th(\mathbb{K}_S; \mathcal{L}_1) = S$ as claimed.

As $S^{\flat}$ is undecidable, $S$ is undecidable. We conclude $\Th(\mathbb{K}_{T_1}|_{\mathcal{L}_1}; \mathcal{L}_1)$ -- namely, $T_1$ as an $\mathcal{L}_1$-theory -- is hereditarily undecidable, as required. \bs

\begin{cor}\label{endcor}
Let $\mathbb{G}$ be the class of nonempty graphs, and let $\mathbb{K}$ be a class of fields. Suppose $\mathbb{G}$ is uniformly interpretable with parameters in $\mathbb{K}$; then $\TTh(\mathbb{K}; \mathcal{L}_r)$ is hereditarily undecidable. 
\end{cor}

\proof
The theory of nonempty graphs is known to be undecidable {(e.g.\ \cite[Corollary 28.5.3]{friedjarden}, which uses a key result of Lavrov \cite[Theorem 3.3.3]{eltt}. In \cite[Theorem III]{church} Church \& Quine note the theory of a binary symmetric predicate is undecidable too). The theory of nontrivial graphs is \textit{hereditarily} undecidable as it is finitely axiomatised (this is \cite[Corollary 3.4.1]{eltt}). \bs

\section{Pseudo-Algebraically Closed Fields}\label{sect31}

A field $K$ is \emph{pseudo-algebraically closed $($PAC$)$} if every geometrically irreducible $K$-variety has a $K$-rational point. The idea behind this section is ultimately the adaptation of the undecidability of the theory of (perfect) PAC fields, due to Cherlin, van den Dries \&\ Macintyre \cite{cvddmac} and independently Ershov \cite{ershov}. The key to this proof was, given an arbitrary nonempty graph $\Gamma$ to construct a (perfect) PAC field interpreting $\Gamma$, as the theory of all such graphs is hereditarily undecidable. This construction was achieved by designing machinery to encode graphs into projective profinite groups, which are the absolute Galois groups of (perfect) PAC fields exactly.

\subsection{Background} Clearly we need apparatus to discuss profinite groups in a first-order setting. To this end (and following the presentation of \cite[\S 5.1]{chatzidakis}) for a profinite group $G$ we consider a structure $S(G)$ whose underlying set is $\bigsqcup_{N \in \mathcal{N}} G/N$, where $\mathcal{N}$ is the family of open normal subgroups of $G$. The elements of $G/N$ are denoted by $gN$, for $g \in G$. $S(G)$ is a structure  in the $\omega$-sorted language $\mathcal{L}_G = \{\leq, C, P\}$, whose sorts are indexed by positive natural numbers $n$ and where $\leq, C$ are binary relations and $P$ is a ternary relation, as follows: the elements of $S(G)$ of sort $n$ are precisely those $gN$ where $N \in \mathcal{N}$ and $[G:N] \leq n$. We say $gN \leq hM$ if and only if $N \subseteq M$, and $C(gN, hM)$ if and only if $N \subseteq M$ \& $gM = hM$. Finally $P(g_1 N_1, g_2 N_2, g_3 N_3)$ if and only if $N_1 = N_2 = N_3$ \& $g_1 g_2 N_1 = g_3 N_1$. 

The $\mathcal{L}_G$-structure $S(G)$ encodes the inverse system $\{(G/N, \pi_{N,M}) \mbox{ : } N, M \in \mathcal{N}\mbox{, } N \subseteq M\}$ precisely and hence determines $G$ uniquely, as $G = \varprojlim G/N$. $S(G)$ is the \textit{complete inverse system associated to $G$}. By \cite[pp.\! 979--980]{chatzidakis}, the class of $\mathcal{L}_G$-structures of the form $S(G)$ for some profinite group $G$ is axiomatisable in the class of all such structures. When $G = G_K = \Gal(K^s/K)$ is the absolute Galois group of a field $K$, this inverse system takes on a new light. Indeed, $S(G_K) = \bigsqcup \Gal(L/K)$ is the union over finite Galois extensions over $K$. The group epimorphisms encoded by $C$ now correspond to the restriction maps $\mbox{res}: \Gal(M/K) \rightarrow \Gal(L/K)$ (for when $K \subseteq L \subseteq M$). If $K/E$ is \textit{regular} (that is to say, $E$ is separably algebraically closed in $K$ and $K/E$ is separable), then the corresponding map of absolute Galois groups $\pi: G_K \twoheadrightarrow G_E$ is surjective, hence $\iota: S(G_E) \hookrightarrow S(G_K)$ by $\iota(N) = \pi^{-1}(N)$, where $N$ is an open normal subgroup of $G_E$.

The language with constants for $Z \subseteq S(G)$ is denoted $\mathcal{L}_G(Z)$. {In this language we will only consider \emph{bounded} formulae and sentences, i.e.\ variables $x$ have a prescribed sort, and quantification over $x$ ranges over the variable's sort. E.g.\ if $x$ is of sort $n$, ``$\exists x$'' is ``there exists $x$ of sort $n$''. The theory $\Th(S(G); \mathcal{L}_G)$ is known as the \emph{cotheory} of $G$. After reconstructing the basic tenants of model theory as ``comodel theory'' for profinite groups in this setting, Cherlin, van den Dries \&\ Macintyre state the following results:

\begin{thm}\label{bitsbobs}
Let $K, L$ be fields with common subfield $E$, such that $K/E$ \&\ $L/E$ are regular. If $K \equiv_E L$ then $S(G_K) \equiv_{S(G_E)} S(G_L)$.

\noindent There is {a \emph{recursive} `translation'} map $-^*: \mbox{\emph{Sent}}(\mathcal{L}_G) \rightarrow \mbox{\emph{Sent}}(\mathcal{L}_r)$ such that if $\varphi$ is a $\mathcal{L}_G$-sentence, then $S(G_K) \models \varphi \iff K \models \varphi^*$, for any field $K$.

\noindent Furthermore, let $\psi(\overline{x})$ be a $\mathcal{L}_G$-formula. Then there is a $\mathcal{L}_r$-formula $\psi^*(\overline{y})$ such that for a tuple $\overline{a}$ of the right sorts in $S(G_K)$, we have
$$(S(G_K), \overline{a}) \models \psi \quad \iff \quad (K, \overline{b}) \models \psi^*,$$
where $\overline{b}$ is a tuple of elements of $K$ encoding $\overline{a} \in S(G_K)$ in a suitable way.
\end{thm}

\proof
{For the former two, see \cite[Theorem 5.9 (1) \& (2)]{chatzidakis} (though \emph{loc.\! cit}.\ does not claim $-^*$ is \emph{recursive}. This is clarified and proven in \cite[Appendix C]{arno}; see \emph{Corollary C.10 ibid}.\ specifically).}

For the latter, {(whose phrasing is taken from \cite[Theorem 6.1.1]{philip}), see \cite[Theorem 5.9 (3)]{chatzidakis} and \cite[Lemma C.8]{arno}}. \bs

\begin{cor}\label{cortocor}
\emph{\cite[Proposition 33]{cvddmac}; \cite[Theorem 5.13]{chatzidakis}.} Let $K_1, K_2$ be PAC fields, separable over a common subfield $E$. Then $K_1 \equiv_E K_2$ if and only if $K_1$ and $K_2$ have the same degree of imperfection, there exists $\theta \in G_E$ such that $\theta(K_1 \cap E^s) = K_2 \cap E^s$, and if $S\Theta : S(G_{K_1 \cap E^s}) \rightarrow S(G_{K_2 \cap E^s})$ is the isomorphism induced by $\theta$, then the partial map $S\Theta: S(G_{K_1}) \rightarrow S(G_{K_2})$ with domain $S(G_{K_1 \cap E^s})$ is $\mathcal{L}_G$-elementary. \bsn
\end{cor}

Suppose $K$ is a PAC field with prime subfield $\mathbb{F}$. Define $K_0 = K \cap \widetilde{\mathbb{F}}$; note $K/K_0$ is regular, hence there is an epimorphism $G_{K} \twoheadrightarrow G_{K_0}$ and thus embedding $S(G_{K_0}) \hookrightarrow S(G_K)$. In the style of Koenigsmann \cite[p.\! 935]{koenis}, denote by $\Th^{alg}(K)$ the subset of $\Th(K; \mathcal{L}_r)$ axiomatised by dictating which monic irreducible one variable polynomials over $\mathbb{F}$ do and do not have a root in $K$ {(see \cite[p.\! 935]{koenis}; note $\Th^{alg}(K) \subseteq \Th(K_0; \mathcal{L}_r)$ but equality does not, in general, hold)}. We have the following corollary to \emph{Corollary \ref{cortocor}}, from utilising the \textit{Compactness Theorem}:

\begin{remark}\label{starrem}
By $\Th(S(G_{K}); \mathcal{L}_G(S(G_{K_0})))^*$ below we mean the set of $\mathcal{L}_r$-sentences obtained from the recursive translation map $-^*$ of \emph{Theorem \ref{bitsbobs}}. Indeed, let $\varphi \in \mbox{Form}(\mathcal{L}_G)$, $\overline{a} \in S(G_{K_0})$, and $S(G_K) \models \varphi(\overline{a})$. Let $\overline{b} \in K_0$ encode $\overline{a} \in S(G_{K_0})$ (by which we mean they are \emph{compatible} in the sense of \cite[Definition C.5]{arno}); then $K \models \varphi^*(\overline{b})$ (this is \cite[Proposition C.9]{arno}) and $\varphi^*(\overline{b}) \in \Th(S(G_{K}); \mathcal{L}_G(S(G_{K_0})))^*$. \sqn
\end{remark}

\begin{cor}\label{mainpaccor}
A PAC field $K$ is axiomatised by the following first-order $(\mathcal{L}_r$-$)$axiom scheme:
\begin{enumerate}
    \item The characteristic and degree of imperfection of $K$;
    \item The PAC field axioms, denoted \emph{\ttfamily PAC};
    \item $\TTh^{alg}(K)$;
    \item $\TTh(S(G_{K}); \mathcal{L}_G(S(G_{K_0})))^*$. \bsn
\end{enumerate}
\end{cor}

A result of van den Dries \& Lubotzky \cite[\S 4.8, Proposition]{lubotzky} describes a correspondence between projective profinite groups and (perfect) PAC fields, via their absolute Galois groups. We would like to use the following improvement:

\begin{thm}\label{genlv}
Let $L/K$ be a Galois extension, $G$ a projective profinite group, and $\alpha: G \rightarrow \Gal(L/K)$ an epimorphism. Then $K$ has an extension $E$ with arbitrary degree of imperfection that is PAC, linearly disjoint from $L$, and there exists an isomorphism $\gamma: G_{E} \rightarrow G$ such that $\alpha \circ \gamma = \mbox{res}_L$.
\end{thm}

{\proof%
\cite[Theorem 23.1.1]{friedjarden}; cf.\! \cite[Proposition 38]{cvddmac}.\bs}

Both \cite{cvddmac} \& \cite{ershov} used the characterisation of the absolute Galois groups of (perfect) PAC fields as projective profinite groups to encode the theory of nonempty graphs into the theory of such PAC field structures. Both proofs have been combined and presented by Fried \&\ Jarden \cite[Chapter 28]{friedjarden} and it is this proof we will reference. Note that, as Ershov remarks, this proof technique demonstrated already the hereditary undecidability of the theory of perfect PAC fields (\cite[p.\! 260]{ershov}). We are aiming for something more general: that every finite subtheory of \emph{any given} PAC field is undecidable.

In {\cite[\S 28.6 -- \S 28.8]{friedjarden}} it is outlined precisely how one can assign a graph $\Gamma_G$ to every profinite group $G$ using two finite groups as parameters (denoted $D$ and $W$): {we present this assignment now. This will be referenced in \S \ref{theprcfields} \& \S \ref{stuffforppc} too.}

\begin{cons}\label{gconstruction}
Given a profinite group $G$, define the graph $\Gamma_G = (A_G, R_G)$ where the set of vertices $A_G$ is the set of open $N \triangleleft G$ such that $G/N \iso D$, and the edge relation $R_G$ is the set of pairs $(N_1, N_2) \in A_G \times A_G$ such that {$N_1 N_2 = G$} and there exists an open $M \triangleleft G$ such that $M \leq N_1 \cap N_2$ and $G/M \iso W$. {Furthermore, there are conditions on the finite groups $D$ and $W = U \rtimes (D \times D)$ that guarantee the surjectivity of the map $G \mapsto \Gamma_G$, known by Fried \& Jarden as \emph{the graph conditions}. They are (from \cite[Definition 28.7.2]{friedjarden}):
\begin{itemize}
    \item[(G1)] $D$, $U$ have no composition factors in common.
    \item[(G2)] For each finite set $I$ and epimorphism $\pi : D^I \rightarrow D$, there exists $i \in I$ such that $\ker(\pi) = \ker(\pi_i)$, where $\pi_i$ is the projection map to the $i$th coordinate of $D^I$.
    \item[(G3)] The intersection of all maximal subgroups of $W$ is trivial, as is the intersection of all maximal subgroups of $D$.
    \item[(G4)] For each embedding $\theta' : D \times D \rightarrow W$ as a semidirect complement (i.e.\ $U \cdot \theta'(D \times D) = W$, and $U \cap \theta'(D \times D) = 1$) and for each nontrivial $N \triangleleft U$, neither factor $D_1 = D \times \{1\}$, $D_2 = \{1\} \times D$ acts trivially {via conjugation on $N$ through $\theta'$.}\footnote{{This is to say there exists $\eta \in N$, $d \in D_1 \cup D_2$ such that $\eta^{\theta'(d)} = \theta'(d) \cdot \eta \cdot \theta'(d)^{-1} \neq \eta$.}}
\end{itemize}
Suppose $D, W$ satisfy these conditions and are in a split short exact sequence:
\[
\begin{tikzcd}
1 \arrow[r] & U \arrow[r] & W \arrow[r, "\lambda", shift left] & D \times D \arrow[l, "\theta", shift left] \arrow[r] & 1.
\end{tikzcd}
\]

{Following the notation of \cite[\S28.8]{friedjarden},} let $\Gamma = (A, R)$ and consider the profinite group $D^A \times W^R$ and the canonical coordinate projections 
\begin{align*}
    &\pi_A : D^A \times W^R \rightarrow D^A; &&\pi_A(\overline{d}, \overline{w}) = \overline{d};\\
    &\pi_a : D^A \times W^R \rightarrow D; &&\pi_a(\overline{d}, \overline{w}) = d_a \quad \mbox{for } a \in A;\\
    &\pi_r : D^A \times W^R \rightarrow {W}; &&\pi_r(\overline{d}, \overline{w}) = w_r \quad \mbox{for } r \in R;    
\end{align*}

Define $G_{\Gamma} \defeq \{(\overline{d}, \overline{w}) \in D^A \times W^R \mbox{ : } r = (a, b) \in R \Rightarrow \lambda(w_r) = (d_a, d_b)\}$. This is a closed subgroup of $D^A \times W^R$, hence is profinite. The graphs $\Gamma,\, \Gamma_{G_{\Gamma}}$ are indeed isomorphic {(\cite[Proposition 28.8.3]{friedjarden})}, and
\[
\begin{tikzcd}
1 \arrow[r] & U^R \arrow[r] & G_{\Gamma} \arrow[r, "\pi_A|_{G_{\Gamma}}"] & D^A \arrow[r] & 1
\end{tikzcd}
\]
is a split short exact sequence {(\cite[Lemmas 28.8.1 \& 28.8.2]{friedjarden}).}} \bs
\end{cons}

For any field $K$, consider specifically the graph $\Gamma_{G_{K}}$: its structure is determined by the absolute Galois group of $K$, which is `seen' in some sense by $K$ through {encoding} finite Galois extensions $K \subset L \subset K^s$ in $K$.

\begin{definition}
{Define the graph $\Gamma_K = (A_K, R_K)$ by}
\begin{align*}
    A_{K} &= \{L \mbox{ : $L/{K}$ is Galois and } \Gal(L/{K}) \iso D\},\\
    R_{K} &= \{(L_1, L_2) \in A_{K} \times A_{K} \mbox{ : } L_1 \cap L_2  = {K} \mbox{ and } \exists N/{K} \mbox{ Galois s.t. }\\
    &\hspace{52mm} L_1L_2 \subseteq N \mbox{ and } \Gal(N/{K}) \iso W\}.
\end{align*}
\end{definition}

\noindent Through the map $L \mapsto G_L$ we see $\Gamma_{K}$ and $\Gamma_{G_{K}}$ are isomorphic. 

\vspace{2mm}
\begin{cons}\label{anothcons}
\cite[pp.\ 694--695]{friedjarden}. For $l \in \Z_{>0}$ let $f_{\overline{x}}(T) = T^l + x_1 T^{l-1} + \dots + x_l$, and if $\overline{a} \in K^l$, denote the splitting field of $f_{\overline{a}}(T)$ over $K$ by $K_{\overline{a}}$. Let $H$ be a finite group,\pagebreak\\ and $\alpha_{l, H}(\overline{x})$ be {an} $\mathcal{L}_r$-formula such that{, for $\overline{a} \in K^l$,}
$$K \models \alpha_{l, H}(\overline{a}) \iff f_{\overline{a}}(T) \mbox{ is separable and } \Gal(K_{\overline{a}}/K) \iso H.$$

Using this and a finite group $E$, with $l = |H|$ {one may construct an $\mathcal{L}_r$-formula $\rho_{H, E}$ such that for $\overline{b}, \overline{c} \in K^l$,}
\begin{align*}
K \models \rho_{H, E}(\overline{b}, \overline{c}) \iff \,& K_{\overline{b}} \cap K_{\overline{c}} = K\mbox{, } K \models \alpha_{l, H}(\overline{b}) \land \alpha_{l, H}(\overline{c})\mbox{, and}\\
&{K \models \exists \overline{z}\,( \alpha_{|E|, E}(\overline{z})\land \mbox{``}K_{\overline{b}} \subseteq K_{\overline{z}}\mbox{''} \land \mbox{``}K_{\overline{c}} \subseteq K_{\overline{z}}\mbox{''}).}    
\end{align*}

We may then define a recursive translation map $-' : \mbox{Form}(\mathcal{L}_{gr}) \rightarrow \mbox{Form}(\mathcal{L}_r)$; $\phi \mapsto \phi'$ by the following rules:
\begin{itemize}
    \item $R(X, Y) \mapsto (R(X, Y))' = \rho_{H, E}(\overline{x}, \overline{y})$;
    \item $\neg \varphi \mapsto \neg(\varphi')$;
    \item $\varphi_1 \land \varphi_2 \mapsto (\varphi_1') \land (\varphi_2')$;
    \item $\exists x (\varphi) \mapsto \exists \overline{x}\, (\alpha_{l, H}(\overline{x}) \land \varphi'(\overline{x}))$. \sq
\end{itemize}
\end{cons}

\begin{lem}\label{citey}
Let $K$ be a field; $\Gamma_K$ is interpretable in $K$.
\end{lem}

\proof
Recalling \emph{Definition \ref{interp}}, set $\mathcal{L}_1 = \mathcal{L}_r$, $\mathcal{L}_0 = \mathcal{L}_{gr}$, $n = |D|$, $\delta(\overline{x}) = \alpha_{n, D}(\overline{x})$, and $f : \alpha_{n, D}(K^n) \rightarrow \Gamma_{K}$; $\overline{a} \mapsto K_{\overline{a}}$. This is indeed surjective, as if $L/K$ is Galois with $\Gal(L/K) \iso D$, it is the splitting field of a degree $n = |D|$ monic separable polynomial $f_{\overline{a}}(T) = T^n + a_1 T^{n-1} + \dots + a_n$ over $K$, and hence $K \models \alpha_{n, D}(\overline{a})$.

Note \emph{Definition \ref{interp} (2)} is satisfied by \emph{Construction \ref{anothcons}} with $E = W$, and condition $(\dagger)$ is confirmed in \cite[p.\ 695]{friedjarden}. \bs

\begin{remark}
Notice the above interpretation is \emph{uniform} in the sense that the $\mathcal{L}_r$-formula $\alpha_{n,D}(\overline{x})$ and the map $-' : \mbox{Form}(\mathcal{L}_{gr}) \rightarrow \mbox{Form}(\mathcal{L}_r)$ of \emph{Construction \ref{anothcons}} do not depend on $K$ or $\Gamma_K$. \sq
\end{remark}

\enlargethispage{\baselineskip}
\noindent We conclude the subsection with some definitions.

\begin{definition}
Let $G$ be a profinite group. The intersection of all maximal open subgroups of $G$ is a normal closed subgroup of $G$ called the \emph{Frattini group of $G$} and denoted $\Phi(G)$.
\end{definition}

\begin{definition}
A homomorphism of profinite groups $\varphi : H \rightarrow G$ is a \emph{Frattini cover} if $\varphi$ is surjective and $\ker(\varphi) \leq \Phi(H)$.
\end{definition}

Given a profinite group $G$, one can partially order the epimorphisms of profinite groups onto $G$: i.e.\ if $\theta_i : H_i \rightarrow G$ for $i=1,2$ are epimorphisms, $\theta_2$ is \emph{larger} than $\theta_1$ if there is an epimorphism $\theta : H_2 \rightarrow H_1$ such that $\theta_1 \circ \theta = \theta_2$. 

\begin{prop}\label{comesafter}
\emph{\cite[Proposition 22.6.1]{friedjarden}.} Every profinite group $G$ has an associated projective group\footnote{Sometimes known as the \emph{Frattini hull of $G$.}} $\widetilde{G}$ and a Frattini cover $\widetilde{\varphi} : \widetilde{G} \rightarrow G$, unique up to isomorphism, called the \emph{universal Frattini cover}, satisfying the following equivalent conditions:
\begin{itemize}
    \item $\widetilde{\varphi}$ is the largest Frattini cover of $G$;
    \item if $\widetilde{G}'$ is a projective profinite group and $\lambda : \widetilde{G}' \twoheadrightarrow G$ is an epimorphism, then $\lambda$ is larger than $\widetilde{\varphi}$. \bs
\end{itemize}
\end{prop}

\subsection{Results} This concludes the pouring of the foundations; we are ready to begin construction. Let $K$ be a PAC field and $\Sigma \subseteq \Th(K;\mathcal{L}_r)$ a finite subtheory. To prove finite undecidability, finding parameters $D_{\Sigma}, U_{\Sigma}, W_{\Sigma}$ such that the Fried-Jarden graph machinery still operates correctly, and does not `interfere' with the part of the absolute Galois group of $K$ axiomatised by $\Sigma$, will be the crucial step. Recall we require $D_{\Sigma}$ and $W_{\Sigma} = U_{\Sigma} \rtimes (D_{\Sigma} \times D_{\Sigma})$ to be finite groups with the following properties:
\begin{itemize}\label{propertiesbegin}
    \item[(G1)] $D_{\Sigma}$, $U_{\Sigma}$ have no composition factors in common.
    \item[(G2)] For each finite set $I$ and epimorphism $\pi : D_{\Sigma}^I \rightarrow D_{\Sigma}$, there exists $i \in I$ such that $\ker(\pi) = \ker(\pi_i)$, where $\pi_i$ is the projection map to the $i$th coordinate of $D_{\Sigma}^I$.
    \item[(G3)] The Frattini subgroups $\Phi(W_{\Sigma}) = \Phi(D_{\Sigma}) = 1$.
    \item[(G4)] For each embedding $\theta' : D_{\Sigma} \times D_{\Sigma} \rightarrow W_{\Sigma}$ as a semidirect complement and for each nontrivial $N \triangleleft U_{\Sigma}$, neither factor $D_1 = D_{\Sigma} \times \{1\}$, $D_2 = \{1\} \times D_{\Sigma}$ acts trivially {via conjugation on $N$ through $\theta'$.}
\end{itemize}

\label{beforegraph}{By the \emph{Compactness Theorem}, there exists a finite set of $\mathcal{L}_r$-sentences $\Delta$ such that $\Delta \models \Sigma$ and $\Delta = \Delta_1 \cup \Delta_2 \cup \Delta_3 \cup \Delta^*_4$, where $\Delta_1$ is a finite subset of $\mathcal{L}_r$-sentences specifying the characteristic and degree of imperfection of $K$ (\emph{Corollary \ref{mainpaccor} (1)}), $\Delta_2$ is a finite subset of {\ttfamily PAC} (\emph{Corollary \ref{mainpaccor} (2)}), $\Delta_3$ is a finite subset of $\Th^{alg}(K)$ (\emph{Corollary \ref{mainpaccor} (3)}), $\Delta_4^*$ is a finite subset of $\Th(S(G_K); \mathcal{L}_G(S(G_{K_0})))^*$ (\emph{Corollary \ref{mainpaccor} (4)}), and $\Delta_4$ is a finite subset of $\Th(S(G_K); \mathcal{L}_G(S(G_{K_0})))$ with $\varphi \in \Delta_4 \Leftrightarrow \varphi^* \in \Delta_4^*$.}

Let $\Lambda$ be the set of universal sentences of $\Delta_3$, and let $\Kbad$ be the join of minimal Galois extensions $F/K$ with $F \models \neg \lambda$ for $\lambda \in \Lambda$. {Let $\overline{a}_{\Sigma} \in S(G_{K_0}) \subset S(G_K)$ be a finite tuple of elements such that $\Delta_4$ is a set of finitely many $\mathcal{L}_G(\overline{a}_{\Sigma})$-sentences. Fix $n_{\Sigma} \in \N$ such that $S_1, \dots, S_{n_{\Sigma}}$ is the smallest consecutive sequence of sorts involving the sentences of $\Delta_4$. (Each sentence $\varphi$ has finitely many occurrences of the symbols $\leq, C, P$, finitely many constant symbols, and finitely many bounded variables. Hence there exists $n_{\varphi} \in \N$ such that the $\mathcal{L}_G(\overline{a}_{\Sigma})$-symbols and variables occurring in $\varphi$ occur in the sorts $S_1, \dots, S_{n_{\varphi}}$.)} Let $\widehat{p}$ be the smallest prime larger than $n_{\Sigma} + |\Gal(\Kbad/K)|$, $P_{\Sigma}$ the set of primes $\{2, 3, \dots, \widehat{p}\}$, and $\mathcal{C}_{\Sigma}$ the formation of finite groups whose order is necessarily a product of powers of primes of $P_{\Sigma}$ (including trivial powers). Define $G_{P_{\Sigma}}$ to be the maximal pro-$\mathcal{C}_{\Sigma}$ quotient of $G_K$.

\begin{cons}\label{graph}
Choose primes $t, r, s$ such that $t > r > s > \widehat{ p}$, $r \equiv 1$ mod $s$, and {$t \equiv 1$ mod $s$} (this can be done by \textit{Dirichlet's Theorem on Arithmetic Progressions}). Let $U_{\Sigma} = C_t$ (the multiplicative cyclic group of order $t$) and $D_{\Sigma} = C_r \rtimes_{\iota} C_s$, where as $r \equiv 1$ mod $s$, there is an embedding $\iota: C_s \hookrightarrow C_{r-1} \iso \Aut(C_r)$ which determines the group operation:
$$(c_1, d_1) \cdot (c_2, d_2) \defeq (c_1 \iota(d_1)(c_2), d_1 d_2).$$

Let $\gamma$ be a generator of $C_r$ and $\beta$ be a generator of $C_s$, and consider these both as elements of $D_{\Sigma}$. (Abusing notation, write $\gamma$ for $(\gamma, 1)$ and $\beta$ for $(1, \beta)$.) Then calculation shows in $C_r \rtimes_{\iota} C_s$ we have the formula
$$\beta^a \cdot \gamma^b = \gamma^{b p(r)^{ak}} \cdot \beta^a,$$

where $r - 1 = ks$ and $p(r)$ is a primitive root modulo $r$ (i.e.\ a generator of $(\Z/r\Z)^*$, the multiplicative group of integers modulo $r$). One can check from this that elements of the form $\gamma^i \cdot \beta^j$ do not commute with one of $\gamma$ or $\beta$ when $0 \leq i < r, 0 \leq j < s$, and $i +j > 0$, hence $D_{\Sigma}$ is centreless. As $t \equiv 1$ mod $s$, there is an embedding $C_s \hookrightarrow C_{t-1} \iso \Aut(C_t)$ which can be extended to a group homomorphism
$$t: D_{\Sigma} = C_r \rtimes C_s \twoheadrightarrow C_s \hookrightarrow \Aut(C_t) = \Aut(U_{\Sigma}).$$

Let $W_{\Sigma} = U_{\Sigma} \rtimes_{\tau} (D_{\Sigma} \times D_{\Sigma})$ under the homomorphism $\tau: D_{\Sigma} \times D_{\Sigma} \rightarrow \Aut(U_{\Sigma})$; $(x, y) \mapsto t(x)t(y)$, and note the following is a split exact sequence:
\begin{equation}
\begin{tikzcd}
1 \arrow[r] & U_{\Sigma} \arrow[r] & W_{\Sigma} \arrow[r, "\lambda", shift left] & D_{\Sigma} \times D_{\Sigma} \arrow[l, "\theta", shift left] \arrow[r] & 1.
\end{tikzcd}\tag*{$\blacksquare$}
\end{equation}
\end{cons}

\begin{lem}\label{dsig}
{With $D_{\Sigma}, U_{\Sigma}, W_{\Sigma}$ from \emph{Construction \ref{graph}}, (G1)--(G4) are satisfied.}
\end{lem}

\proof
We follow \cite[Example 28.7.4]{friedjarden} as much as possible.

(G1) is satisfied: the composition factors of $D_{\Sigma}$ are $C_r$ and $C_s$, distinct to the (unique) composition factor $C_t$.

(G2) is satisfied: let $I$ be a finite set and $\pi : D_{\Sigma}^I = \prod_{i \in I} (D_{\Sigma})_i \rightarrow D_{\Sigma}$ an epimorphism. We may suppose $\ker(\pi) \cap (D_{\Sigma})_i$ is a \underline{proper} (normal) subgroup of $D_{\Sigma}$ for each $i \in I$; otherwise if $(D_{\Sigma})_j \leq \ker(\pi)$ for some $j \in I$ we may consider the epimorphism $\pi' : D_{\Sigma}^{I \setminus \{j\}} \rightarrow D_{\Sigma}$ and apply induction to find $i \in I\setminus \{j\}$ such that $\ker(\pi) = \ker(\pi_i)$.

For each $i \in I$, $\ker(\pi) \cap (D_{\Sigma})_i$ is thus either 1 or a cyclic subgroup of $(D_{\Sigma})_i$ of prime order. As $\ker(\pi) \cap (D_{\Sigma})_i$ is normal, then $\ker(\pi) \cap (D_{\Sigma})_i = 1$ or $\langle \gamma_i\rangle$ where $\gamma_i \in (D_{\Sigma})_i$ is of order $r$. If $\ker(\pi) \cap (D_{\Sigma})_i = \langle \gamma_i\rangle$ for all $i \in I$, then by the \textit{First Isomorphism Theorem for Groups} the centre of $D_{\Sigma}$ is nontrivial -- a contradiction. WLOG suppose $\ker(\pi) \cap (D_{\Sigma})_1 = 1$. By the \textit{First Isomorphism Theorem for Groups}, $\pi((D_{\Sigma})_1) = D_{\Sigma}$, and every element of $\pi((D_{\Sigma})_i)$ for $i \neq 1$ commutes with every element of $D_{\Sigma}$. This is a contradiction unless $I = \{1\}$ (as we have assumed $\ker(\pi) \cap (D_{\Sigma})_i$ is a \textit{proper} subgroup of $D_{\Sigma}$ for each $i \in I$).

(G3) is satisfied: from the formula $\beta^a \cdot \gamma^b = \gamma^{b p(r)^{ak}} \cdot \beta^a$ one can calculate
$$(\beta \cdot \gamma)^n = \gamma^{\sum_{i=1}^n p(r)^{ik}} \cdot \beta^n.$$

The order of $\beta\gamma$ is thus $s$, as $s$ is the smallest power $n$ for which $\beta^n = 1$, and recalling $sk = r-1$, $p(r)^{r-1} \equiv 1$ mod $r$, note:
$$\sum_{i=1}^s p(r)^{ik} \equiv \sum_{i=0}^{s-1} (p(r)^k)^i \equiv \frac{(p(r)^k)^s - 1}{p(r)^k - 1} \equiv 0 \mbox{ mod } r.$$

Hence $\langle \beta\gamma\rangle$ and $\langle \beta \rangle$ both have index $r$ in $D_{\Sigma}$, and $\langle \beta \gamma \rangle \cap \langle \beta \rangle = 1$. By definition $\Phi(D_{\Sigma}) \leq \langle \beta \gamma \rangle \cap \langle \beta \rangle$, hence $\Phi(D_{\Sigma}) = 1$ and $\Phi(W_{\Sigma}) \leq U_{\Sigma}$. However $\theta(D_{\Sigma} \times D_{\Sigma})$ has index $t$ in $W_{\Sigma}$, thus $\Phi(W_{\Sigma}) \leq U_{\Sigma} \cap \theta(D_{\Sigma} \times D_{\Sigma}) = 1$.

Finally, (G4) is satisfied: let $\theta'$ be an embedding of $D_{\Sigma} \times D_{\Sigma}$ into $W_{\Sigma}$ as a semidirect complement. Since the orders of $U_{\Sigma}$ and $D_{\Sigma} \times D_{\Sigma}$ are relatively prime, $\theta'(D_{\Sigma} \times D_{\Sigma})$ is conjugate to $\theta(D_{\Sigma} \times D_{\Sigma})$ by the \textit{Schur-Zassenhaus Lemma}\footnote{{The formulation we use is \cite[Lemma 22.10.1]{friedjarden}.}}. By (G2), $D_{\Sigma} \times D_{\Sigma}$ has a unique factorisation as a direct product of two copies of $D_{\Sigma}$, meaning we may canonically write $D_{\Sigma} \times D_{\Sigma} = D_1 \times D_2$. We have, for $w \in W_{\Sigma}$, $v \in U_{\Sigma}$ (also denoted $v$ as an element of $W_{\Sigma}$), and $d \in D_1 \cup D_2$:
\begin{align*}
    v^{\theta'(d)} &= v^{\theta(d)^w} \\
    &= \theta(d)^w \cdot v \cdot (\theta(d)^w)^{-1} = (w \cdot \theta(d) \cdot w^{-1}) \cdot v \cdot (w \cdot \theta(d) \cdot w^{-1})^{-1}\\
    &= (w \cdot \theta(d) \cdot w^{-1}) \cdot v \cdot (w \cdot \theta(d)^{-1} \cdot w^{-1}) = w \cdot [\theta(d) \cdot [w^{-1} \cdot v \cdot w] \cdot \theta(d)^{-1}] \cdot w^{-1}\\
    &=w \cdot \tau(d)(w^{-1} \cdot v \cdot w) \cdot w^{-1}.
\end{align*}

Note as $U_{\Sigma}$ is normal, $\overline{v} \defeq w^{-1} \cdot v \cdot w \in U_{\Sigma}$. If $w \cdot \tau(d)(w^{-1} \cdot v \cdot w) \cdot w^{-1} = v$, then $\tau(d)(\overline{v}) = \overline{v}$. However we may choose $v \in U_{\Sigma}$ and $d \in D_1 \cup D_2$ such that $\overline{v} \neq 1$ and $\tau(d) \in \Aut(U_{\Sigma})$ has $\tau(d)(\overline{v}) \neq \overline{v}$. Hence neither of $D_1$ or $D_2$ act trivially via {conjugation  on $U_{\Sigma}$ through $\theta'$}, as desired. \bs

Given a graph $\Gamma = (A, R)$, {consider the profinite group $G_{\Gamma}$ from \emph{Construction \ref{gconstruction}} using $D_{\Sigma}$ and $W_{\Sigma}$} (we suppress the $\Sigma$ notation in $G_{\Gamma}$). There is a split exact sequence
\[
\begin{tikzcd}
1 \arrow[r] & U^R_{\Sigma} \arrow[r] & G_{\Gamma} \arrow[r] & D^A_{\Sigma} \arrow[r] & 1.
\end{tikzcd}
\]

Consider the cotheory of $G_{\Gamma}$: by design, there are no proper open normal subgroups $N \triangleleft G_{\Gamma}$ of index less than or equal to $\widehat{p}$. Hence the only element of $S(G_{\Gamma})$ of sort $n \leq \widehat{p}$ is $1 G_{\Gamma}$. Note also, $\Gamma \iso \Gamma_{G_{\Gamma}}$ from \emph{Construction \ref{gconstruction}.} \label{prof}Let $\widetilde{G_{\Gamma}}$ be the universal Frattini cover of $G_{\Gamma}$ with map $f : \widetilde{G_{\Gamma}} \rightarrow G_{\Gamma}$. By \cite[Corollary 54]{cvddmac} the sorts $S_1, \dots, S_{\widehat{p}}$ of $S(\widetilde{G_{\Gamma}})$ remain trivial, and $\Gamma_{\widetilde{G_{\Gamma}}} \iso \Gamma_{G_{\Gamma}}$ by \cite[Lemma 28.6.1]{friedjarden}. Let $G^*_{\Gamma} = G_{P_{\Sigma}} \star \widetilde{G_{\Gamma}}$, the free product\footnote{The coproduct in the category of profinite groups. {In \cite[Lemma 22.4.9]{friedjarden} it is proven this is the profinite completion of the free product $\mathbb{G}$ of abstract groups $G_{P_{\Sigma}}, \widetilde{G_{\Gamma}}$, with respect to the collection of $N \triangleleft \mathbb{G}$ of finite index, such that $N \cap G_{P_{\Sigma}}$ is open in $G_{P_{\Sigma}}$ and $N \cap \widetilde{G_{\Gamma}}$ is open in $\widetilde{G_{\Gamma}}$.}} of $G_{P_{\Sigma}}$ and $\widetilde{G_{\Gamma}}$. For reference later we prove the following:

\begin{lem}\label{addedsea}
$S(G_{P_{\Sigma}}), S(G^*_{\Gamma}) \models \Delta_4$. 
\end{lem}

\proof
Consider $S(G_{P_{\Sigma}})$ as $\mathcal{L}_G$-structure: it may be expanded to an $\mathcal{L}_G(\overline{a}_{\Sigma})$-structure. Indeed, for $a_i \in\overline{a}_{\Sigma}$, let $g_i N_i$ be the interpretation of $a_i$ in $S(G_K)$, where $g_i N_i \in G_K/N_i$ is of sort $n_i$, $N_i \triangleleft S(G_K)$. By design, $n_i \leq \widehat{p}$ and there exists $N'_i \triangleleft G_{P_{\Sigma}}$ with $\psi_i : G_K/N_i \xrightarrow{\iso} G_{P_{\Sigma}}/N'_i$ by the \emph{Isomorphism Theorems for Compact Groups} \cite[p.\ 5]{friedjarden}. Define the interpretation of $a_i$ in $S(G_{P_{\Sigma}})$ to be $\psi_i(g_i N_i)$. 

In fact, there is a correspondence between closed $N_k \triangleleft G_K$ of index $k \leq \widehat{p}$, and closed $N'_k \triangleleft G_{P_{\Sigma}}$ of index $k \leq \widehat{p}$, with $G_K/N_k \iso G_{P_{\Sigma}}/N'_k$, by \cite[Theorem 2.28(ii)]{rotman} (note the quotient map $G_K \twoheadrightarrow G_{P_{\Sigma}}$ is continuous, so Rotman's proof holds for \emph{closed} subgroups). Therefore if $\varphi$ is a sentence of $\Delta_4$, $S(G_{P_{\Sigma}}) \models \varphi$, as $S(G_K) \models \varphi$ and the $\mathcal{L}_G(\overline{a}_{\Sigma})$-symbols and variables of $\varphi$ occur in the sorts $S_1, \dots, S_{\widehat{p}}$. 

Similarly, $S(G^*_{\Gamma}) \models \Delta_4$: if $N_i \triangleleft G^*_{\Gamma}$ is of index $\leq \widehat{p}$, under the projection $G^*_{\Gamma} \twoheadrightarrow G^*_{\Gamma}/N_i$ we must have $\widetilde{G_{\Gamma}} \leq N_i$, as otherwise $\widetilde{G_{\Gamma}}$ would have a proper normal subgroup of index $\leq \widehat{p}$, a contradiction to the above. By \cite[Theorem 2.28(ii)]{rotman} there is a correspondence between closed $N_k \triangleleft G^*_{\Gamma}$ of index $k \leq \widehat{p}$, and closed $N'_k \triangleleft G_{P_{\Sigma}}$ of index $k \leq \widehat{p}$, with $G^*_{\Gamma}/N_k \iso G_{P_{\Sigma}}/N'_k$, (again note the quotient map $G^*_{\Gamma} \twoheadrightarrow \widetilde{G_{\Gamma}}/\langle \widetilde{G_{\Gamma}}\rangle \iso G_{P_{\Sigma}}$ is continuous, so Rotman's proof holds for {closed} subgroups). Hence for any $\mathcal{L}_G(\overline{a}_{\Sigma})$-sentence $\varphi$ with variables over the sorts $S_1, \dots, S_{\widehat{p}}$, $S(G_{P_{\Sigma}}) \models \varphi \Leftrightarrow S(G^*_{\Gamma}) \models \varphi$. We conclude $S(G^*_{\Gamma}) \models \Delta_4$, as desired. \bs

Notice $G^*_{\Gamma}$ is projective: as the maximal pro-$\mathcal{C}_{\Sigma}$ quotient of the projective profinite group $G_K$, $G_{P_{\Sigma}}$ is projective (\cite[Proposition 22.4.8]{friedjarden}), and the free product of projective profinite groups is projective (\cite[Proposition 22.4.10]{friedjarden}). We claim $\Gamma_{G^*_{\Gamma}} \iso \Gamma_{\widetilde{G_{\Gamma}}}$, or more generally:

\begin{lem}\label{graphimp}
Let $H$ be a pro-$\mathcal{C}_{\Sigma}$ group; then $\Gamma_{H \,\star\,\widetilde{G_{\Gamma}}} \iso \Gamma_{\widetilde{G_{\Gamma}}}$. 
\end{lem}

\proof
{A minor adaptation of} \cite[Lemma 28.6.1]{friedjarden}. Consider the quotient map $\pi :H \star\widetilde{G_{\Gamma}} \rightarrow \widetilde{G_{\Gamma}}$, with kernel $\langle H\rangle$ the least normal subgroup of $H \star\widetilde{G_{\Gamma}}$ containing $H$. For open normal subgroups $N \triangleleft \widetilde{G_{\Gamma}}$, there is an isomorphism
$$\widetilde{G_{\Gamma}}/N \iso (H  \star \widetilde{G_{\Gamma}})/\pi^{-1}(N).$$
This yields an embedding $\Gamma_{\widetilde{G_{\Gamma}}} \hookrightarrow \Gamma_{H \,\star\,\widetilde{G_{\Gamma}}}$ along $\pi^{-1}$.

Conversely, if $M \triangleleft H \star \widetilde{G_{\Gamma}}$ is open, and $H  \star \widetilde{G_{\Gamma}}/M \iso D_{\Sigma}$ (resp.\ $W_{\Sigma}$), then consider $\rho : H \star \widetilde{G_{\Gamma}} \twoheadrightarrow D_{\Sigma}$ (resp.\ $W_{\Sigma}$) and {notice $\rho(H) \leq D_{\Sigma}$ (resp.\ $W_{\Sigma}$). Since any quotient of $H$ by an open normal subgroup of $H$ is a $P_{\Sigma}$-group, there exists elements of $\rho(H)$ of the wrong order (by \emph{Lagrange's Theorem}) unless $H \leq \ker(\rho)$. As $\ker(\rho)$ is a normal subgroup of $H \star \widetilde{G_{\Gamma}}$, $\ker(\pi) = \langle H \rangle \leq \ker(\rho) = M$. Therefore $M = \pi^{-1}(\pi(M))$ and thus} $\pi^{-1}$ induces an isomorphism of graphs $\Gamma_{H \,\star\,\widetilde{G_{\Gamma}}} \iso \Gamma_{\widetilde{G_{\Gamma}}}$ as required. \bs

\noindent We are ready to prove: 

\begin{thm}\label{thisolthm}
\textit{For any nontrivial graph $\Gamma$ there exists a PAC field $K_{\Gamma} \supseteq K$ such that $K_{\Gamma} \models \Sigma$ and $\Gamma_{K_{\Gamma}} \iso \Gamma$.}
\end{thm}

\proof
Given $\Gamma$, $\Sigma$, construct the projective profinite group $G^*_{\Gamma}$. Note that $\Gal(\Kbad/K)$ is a quotient of $G_{P_{\Sigma}}$, which is a quotient of $G^*_{\Gamma}$, hence there is an epimorphism $G^*_{\Gamma} \twoheadrightarrow \Gal(\Kbad/K)$. By \emph{Theorem \ref{genlv}} there is a PAC field $K_{\Gamma}$ extending $K$, of the same degree of imperfection as $K$, such that $K_{\Gamma} \cap \Kbad = K$ and $G_{K_{\Gamma}} \iso G^*_{\Gamma}$. Therefore $K_{\Gamma} \models \Delta_1 \cup \Delta_2 \cup \Delta_3 \cup \Delta^*_4$ by construction, hence $K_{\Gamma} \models \Sigma$, and $\Gamma \iso \Gamma_{G_{K_{\Gamma}}} \!\iso \Gamma_{K_{\Gamma}}$ as required. \bs

\begin{cor}\label{mainpacfieldundeccor}
Every PAC field is finitely undecidable.
\end{cor}

\proof
Fix $K$ a PAC field and $\Sigma \subseteq \Th(K; \mathcal{L}_r)$ a finite subtheory. Let $\Gamma$ be a nonempty graph; by \emph{Lemma \ref{citey} \& Theorem \ref{thisolthm}} there exists a PAC field $K_{\Gamma} \models \Sigma$ and $\Gamma$ is interpretable in $K_{\Gamma}$. Furthermore, the interpretation of \emph{Lemma \ref{citey}} is uniform in the sense of \emph{Definition \ref{uniformparam}}; the class of nontrivial graphs is uniformly interpretable in the class of PAC fields satisfying $\Sigma$. We conclude {\ttfamily PAC}$\,\cup\, \Sigma$ is hereditarily undecidable by \emph{Corollary \ref{endcor}}, hence $\Sigma$ is undecidable as required. \bs

\begin{exmp}
Let $\mathbb{K}$ be a countable $\omega$-free PAC field of characteristic $0$ containing $\widetilde{\Q}$. By standard results (cf.\ {\cite[Proposition 9]{koenis}}) $\Th(\mathbb{K};\mathcal{L}_r)$ is decidable, however by \emph{Corollary \ref{mainpacfieldundeccor}}, $\Th(\mathbb{K};\mathcal{L}_r)$ is finitely undecidable. \sq
\end{exmp}

Although it follows from \emph{Corollary \ref{mainpacfieldundeccor}}, with much simpler tools we can prove the following, outlined to the author by E.\ Hrushovski:

\begin{thm}\label{pacfiniteaxiom}
No PAC field is finitely axiomatisable (even among the class of PAC fields of the same characteristic and degree of imperfection).
\end{thm}

(This is to say that, given a PAC field $L$, there does not exist an $\mathcal{L}_r$-sentence $\gamma$ such that for all PAC fields $F$ of the same characteristic and degree of imperfection as $L$, $F \models \gamma \iff F \equiv_{\mathcal{L}_r} L$.)\\

\proof
Assume for the purpose of contradiction there exists a finitely axiomatisable PAC field $K$; in particular, by \emph{Corollary \ref{mainpaccor}} to characterise $K$ among all PAC fields of the same characteristic and degree of imperfection, we need only finite many axioms $\Sigma_1 \subset \Th^{alg}(K)$, $\Sigma_2^* \subset \Th(S(G_K); \mathcal{L}_G(S(G_{K_0})))^*$ (where $\Sigma_2 \subset \Th(S(G_K); \mathcal{L}_G(S(G_{K_0})))$ is finite).

Let $\Lambda$ be the set of universal sentences of $\Sigma_1$; $\Lambda$ specifies the monic irreducible univariate polynomials over $\mathbb{F}$ (the prime subfield) in $\Sigma_1$ which \emph{do not} have a root in $K$. Let $\Kbad$ be the join of minimal Galois extensions $F/K$ with $F \models \neg \lambda$ for $\lambda \in \Lambda$. Then $\Kbad/K$ is finite, $\Kbad \models \neg \Lambda$, and $\Gal(\Kbad/K)$ is a quotient of $G_K$. As there is an epimorphism $G_K \twoheadrightarrow \Gal(\Kbad/K)$, there is an embedding of $\mathcal{L}_G$-structures $S(\Gal(\Kbad/K)) \hookrightarrow S(G_K)$.

{Let $\overline{a} \in S(G_{K_0}) \subset S(G_K)$ be a finite tuple of elements such that $\Sigma_2$ is a set of  (finitely many) $\mathcal{L}_G(\overline{a})$-sentences. Fix $n_k \in \N$ such that $S_1, \dots, S_{n_k}$ is the smallest consecutive sequence of sorts involving the sentences of $\Sigma_2$. (Each sentence $\varphi$ has finitely many occurrences of the symbols $\leq, C, P$, finitely many constant symbols, and finitely many bounded variables. Hence there exists $n_{\varphi} \in \N$ such that the $\mathcal{L}_G(\overline{a})$-symbols and variables occurring in $\varphi$ occur at the sorts $S_1, \dots, S_{n_{\varphi}}$.)} Let $\widehat{p}$ be the smallest prime larger than $n_k + |\Gal(\Kbad/K)|$, $P$ the set of primes $\{2, 3, \dots, \widehat{p}\}$, and $\mathcal{C}$ the formation of finite groups whose order is necessarily a product of powers of primes of $P$ (including trivial powers). As this formation is full, and $G_K$ is projective, the maximal pro-$\mathcal{C}$ quotient of $G_K$ (denoted $G_P$) is also projective (\cite[Proposition 22.4.8]{friedjarden}). Moreover, by construction $\Gal(\Kbad/K)$ is up to isomorphism a quotient of $G_P$, hence there is an epimorphism $G_P \twoheadrightarrow \Gal(\Kbad/K)$.

Let $\mathbb{F}_1(q)$ be the free pro-$q$ group on one generator, where $q$ is a prime larger than $\widehat{p}$. Finally, let $G_P \star \mathbb{F}_1(q)$ be the free product of $G_P$ and $\mathbb{F}_1(q)$. By \cite[Proposition 22.4.10]{friedjarden}, $G_P \star \mathbb{F}_1(q)$ is a projective profinite group. The $\mathcal{L}_G$-theories of $G_P \star \mathbb{F}_1(q)$ and $G_P$ clearly differ; there is an epimorphism from $G_P \star \mathbb{F}_1(q)$ to $\Z / q\Z$, however no such epimorphism from $G_P$ exists by construction. {However $S(G_P), S(G_P \star \mathbb{F}_1(q)) \models \Sigma_2$ by the argument of \emph{Lemma \ref{addedsea}} (replacing $\widetilde{G_{\Gamma}}$ there by $\mathbb{F}_1(q)$).

By \emph{Theorem \ref{genlv}} there exist PAC fields $F_1, F_2 \supseteq K$ such that $G_{F_1} \iso G_P$, $G_{F_2} \iso G_P \star \mathbb{F}_1(q)$, and $F_1 \cap \Kbad = F_2 \cap \Kbad = K$. We conclude that $F_1, F_2 \models \Sigma_1 \cup \Sigma^*_2$, however $F_1 \not\equiv F_2$ as $S(G_{F_1}) \not\equiv_{\mathcal{L}_G} S(G_{F_2})$; a contradiction as required. \bs

As an undecidability result, \emph{Corollary \ref{mainpacfieldundeccor}} is interesting in its own right. With \emph{Theorem \ref{pacfiniteaxiom} \&\ Corollary \ref{mainpacfieldundeccor}} there are connections back to \emph{Problems \ref{1} \& \ref{3}}, modulo a classification-theoretic conjecture:

\begin{customconj}{(Simple Fields)}
\textit{Every infinite simple field is PAC.} \bsn
\end{customconj}

\begin{cor}
Assume the Simple Fields Conjecture. Then every infinite simple field is not finitely axiomatisable, and furthermore is finitely undecidable. \bsn
\end{cor}

\bigskip 
\section{Pseudo-Real Closed Fields}\label{theprcfields}
After considering PAC fields, the next natural step to take is to the \emph{pseudo-real closed} fields; the PAC-analogue of an ordered field. PRC fields were given their modern formulation by Prestel \cite[Theorem 1.2]{prestel}:

\begin{definition}
{A field $K$ is \emph{pseudo-real closed} (PRC) if every geometrically irreducible variety defined over $K$, which has a smooth $\overline{K}$-rational point in each real closure $\overline{K}$ of $K$, has a $K$-rational point.}
\end{definition}

The theory of formally real PRC fields is also (hereditarily) undecidable by the work of Haran. It is Haran's proof we will adapt to determine finite undecidability for all PRC fields. Familiarity with the ideas and notation of \cite{haranjarden} is assumed in this section.

\begin{thm}
\emph{\cite[Theorem 3.1]{haran}.} Let $\Xi$ be a nonempty family of Boolean spaces, and $\mbox{\emph{\ttfamily PRC}}(\Xi)$ the elementary theory of the class of PRC fields $K$ such that $X(K) \in \Xi$. Then $\mbox{\emph{\ttfamily PRC}}(\Xi)$ is (hereditarily) undecidable. \bs
\end{thm}

\subsection{Background} Recall from \cite{haranjarden} the category of \textit{Artin-Schreier structures}; this is the `right' category to consider the Galois theory of PRC fields in, as evidenced by the main result of \cite{haranjarden}:

\begin{thm}\label{bigharanjarden}
\emph{\cite[Theorem 10.4]{haranjarden}.} If $K$ is a PRC field, then $\mathfrak{G}_K$ -- \emph{the absolute Artin-Schreier structure of $K$} -- is projective. Conversely, if $\mathfrak{G}$ is a projective Artin-Schreier structure, there exists a PRC field $K$ such that $\mathfrak{G} \iso \mathfrak{G}_K$. \bsn
\end{thm}

(``Projective'' in the sense of \cite[Lemma 7.5]{haranjarden}, though one can prove this notion coincides with the category-theoretic notion of ``projective'' in the category of Artin-Schreier structures.) Haran \&\ Jarden actually prove a theorem slightly more precise than this. They prove a corresponding \emph{Theorem \ref{genlv}:} given a PRC field $K$, they produce a PRC field extension $E$ whose algebraic part is governable relative to $K$, yet has an almost arbitrary absolute Artin-Schreier structure. 

\begin{thm}\label{genlv2}
\emph{\cite[Theorem 10.2]{haranjarden}.} Let $\mathfrak{G}$ be a projective Artin-Schreier structure. Let $L/K$ be a Galois extension such that $\sqrt{-1} \in L$ and let $\pi : \mathfrak{G} \rightarrow \mathfrak{Gal}(L/K)$ be an epimorphism. Then there exists a PRC extension $E/K$ such that $\mathfrak{G} \iso \mathfrak{G}_E$, and
\[
\begin{tikzcd}
\mathfrak{G} \arrow[r, "\iso"] \arrow[rd, "\pi", swap] & \mathfrak{G}_E \arrow[d, "\mbox{Res}_L"]\\
& \mathfrak{Gal}(L/K)
\end{tikzcd}
\]
is a commutative diagram. \bsn
\end{thm}

\noindent (Recall \cite[Example 3.2]{haranjarden} for the notation $\mathfrak{Gal}(L/K)$.) Immediately, we obtain:

\begin{cor}\label{214}
Let $L/K$ be a Galois extension, $\mathfrak{G}$ a projective Artin-Schreier structure, and $\alpha : \mathfrak{G} \rightarrow \mathfrak{Gal}(L(\sqrt{-1})/K)$ an epimorphism. Then $K$ has an extension $E$ which is PRC, linearly disjoint from $L(\sqrt{-1})$, and there exists an isomorphism $\gamma: \mathfrak{G}_E \rightarrow \mathfrak{G}$ such that $\alpha \circ \gamma = \mbox{Res}_L$. \bsn
\end{cor}

There is a rich theory of (absolute) Artin-Schreier structures one can develop alongside the theory of (absolute) Galois groups, as done by Haran \& Jarden \cite{haranjarden} and Haran \cite{haran, harancohom, harancohom2}.

\begin{thm}\label{thmforprcel}
\emph{{(cf.\ \cite[Theorem 2]{ershovtwotheorem}, \cite[Proposition 3]{jarden}.)}} Let $K_1, K_2$ be PRC fields, {separable} over a common subfield $E$. Then $K_1 \equiv_E K_2$ if and only if $K_1$ \&\ $K_2$ have {the same degree of imperfection}, there exists $\theta \in G_E$ such that $\theta(K_1 \cap E^s) = K_2 \cap E^s$, and if ${S\Theta : S(G_{K_1 \cap E^s}) \rightarrow S(G_{K_2 \cap E^s})}$ is the isomorphism induced by $\theta$, then the partial map $S\Theta : S(G_{K_1}) \rightarrow S(G_{K_2})$ with domain $S(G_{K_1 \cap E^s})$ is $\mathcal{L}_G$-elementary.
\end{thm}

\proofsketch
{The following is Chatzidakis \cite[Theorem 5.13]{chatzidakis}, mutatis mutandis}. Assuming $K_1 \equiv_E K_2$, it is clear there exists $\theta \in G_E$ with $\theta(K_1 \cap E^s) = \theta(K_2 \cap E^s)$, and the partial map $S\Theta : S(G_{K_1}) \rightarrow S(G_{K_2})$ induced by $\theta$ is $\mathcal{L}_G$-elementary on $S(G_{K_1 \cap E^s})$.

Conversely, assume we have the maps $\theta$, $S\Theta$. Moving ${K_1}^s$ by an automorphism extending $\theta$, we may assume WLOG $K_1 \cap E^s = K_2 \cap E^s$, then replace $E$ by $K_1 \cap E^s$ and $\theta$ by the identity map. Note now $K_1, K_2$ are regular extensions of $E$, and $S\Theta$ is the identity on $S(G_E)$. By an analogue of the \textit{Keisler-Shelah Theorem}, there is a nonprincipal ultrafilter $\mathcal{U}$ on an index set $I$ such that $S(G_{K_1^{\mathcal{U}}}) \iso_{S(G_E)} S(G_{K_2^{\mathcal{U}}})$ {(see \cite[(5.12)]{chatzidakis} with \emph{Remarks (2)}, \S5.5 \emph{ibid.})}. Dualising, there is a group homeomorphism $\varphi : G_{K_1^{\mathcal{U}}} \rightarrow G_{K_2^{\mathcal{U}}}$ such that for every $\sigma \in G_{K_1^{\mathcal{U}}}$, $\varphi(\sigma)|_{E^s} = \sigma|_{E^s}$. By Cherlin, van den Dries \&\ Macintyre \cite[\S 3]{cvddmac}\footnote{{Presented also in Fried \& Jarden \cite[Theorem 20.3.3]{friedjarden}.}} (for nonformally real PRC fields) and Jarden \cite[Proposition 3]{jarden} (for formally real PRC fields; cf.\ Ershov \cite{ershovtwotheorem}), this forces $K_1^{\mathcal{U}} \equiv_E K_2^{\mathcal{U}}$, hence $K_1 \equiv_E K_2$ as required. \bs

Note this is elementary equivalence in the language of rings (over a common subfield); \emph{not} the language of \emph{ordered} rings. 

\begin{cor}\label{prcax}
A PRC field $K$ is axiomatised by the following first-order $(\mathcal{L}_r$-$)$axiom scheme:
\begin{enumerate}
    \item The characteristic and degree of imperfection of $K$;
    \item The PRC field axioms, $\mbox{\emph{\ttfamily PRC}}$;
    \item $\TTh^{alg}(K)$;
    \item $\TTh(S(G_K); \mathcal{L}_G(S(G_{K_0})))^*$, where $K_0 = K \cap \mathbb{F}^s$ and $\mathbb{F}$ is the prime subfield of $K$. \phantom{boo}\bsn
\end{enumerate}
\end{cor}

\begin{remark}\label{420}
We emphasise here that if $K$ is a \textit{formally real} PRC field (distinct to the nonformally real PRC fields; the PAC fields) then necessarily $K$ has characteristic $0$. In addition, the orders on $K$ can in a way be `seen' by $\Th(S(G_K); \mathcal{L}_G(S(G_{K_0})))$. Indeed, the space of orders $X(K) = X(K(\sqrt{-1})/K) \iso X(K^s/K)/G_{K(\sqrt{-1})}$, which is homeomorphic to $\mbox{Inv}(G_{K})/G_{K(\sqrt{-1})}$ by \cite[Theorem 2.1]{harancohom2}. \sq
\end{remark}

Now that we have a solid description of the first-order theory of any given PRC field, we are almost ready to conclude its finite undecidability -- what remains are category-theoretic tools, such as Haran's method of transferring the graph constructions in the category of profinite groups to the category of Artin-Schreier structures (in \cite[\S 2]{haran}).

{The following recalls \cite[p.\ 102]{haran}:} fix a Boolean space $X$. For a profinite group $G$, let $\Z\! / 2\Z \times G$ act on the Boolean space $X \times G$ by:\label{73}
$$(x, g)^{(1,h)} = (x, g)^{(\varepsilon, h)} = (x, gh), \qquad \mbox{for } x \in X \mbox{ and } g, h \in G,$$

where $\Z\!/2\Z = \langle \varepsilon\rangle$. Define $d : X \times G \rightarrow \Z\!/2\Z \times G$ by $d(x, g) = \varepsilon$. We have constructed an Artin-Schreier structure
$$F(X, G) = \langle \Z\!/2\Z \times G\mbox{, } G\mbox{, } d : X \times G \rightarrow \Z\!/2\Z \times G\rangle.$$

This describes a (faithful) functor $F(X, -)$ from the category of profinite groups to the category of Artin-Schreier structures, where we extend $F$ to morphisms in the obvious way, with the additional property that the `orbit space' $X(F(X, G))/G \iso X$. Moreover if $\varphi$ is an epimorphism of profinite groups, then $F(X, \varphi)$ is a cover of Artin-Schreier structures. 

\begin{remark}\label{snaky}
\cite[Remark 2.2]{haran}. Let $L/K$ be a Galois extension, $G$ a profinite group, and $X$ a Boolean space; then $\mathfrak{Gal}(L/K) \iso F(X, G)$ if and only if $\sqrt{-1} \in L \setminus K$, $X(K) \iso X$, and there exists a \textit{totally real}\footnote{Each ordering $P \in X(K)$ extends to one on $L_0$.} Galois extension $K \subseteq L_0 \subseteq L$ such that $\mbox{Gal}(L_0/K) \iso G$ and $L = L_0(\sqrt{-1})$. \sq
\end{remark}

We will also use the forgetful functor from the category of Artin-Schreier structures to the category of profinite groups: 
(morphisms are handled in the obvious way)
$$\mbox{Ft}(\mathfrak{G}) = \mbox{Ft}(\langle G, G', X(\mathfrak{G}) \xrightarrow{d} G\rangle) = G.$$

Given an Artin-Schreier structure $\mathfrak{G} = \langle G, G', X(\mathfrak{G}) \xrightarrow{d} G\rangle$, we may perform an analogous graph construction for a fixed Boolean space $X$, and as before the groups $D, W$ mentioned here satisfy conditions (G1)--(G4) of \emph{Construction \ref{gconstruction}}. Let $A^X_{\mathfrak{G}}$ be the set of open normal subgroups of $G$ such that $N \leq G'$, $\mathfrak{G}/N \iso F(X, D)$, and define the binary relation $R^X_{\mathfrak{G}}$ on $A^X_{\mathfrak{G}}$ to be the following set:
$$\{(N_1, N_2) \in A^X_{\mathfrak{G}} \times A^X_{\mathfrak{G}} \mbox{ : } N_1 \!\neq\! N_2\mbox{ \&\ } \exists M \triangleleft G \mbox{ open s.t.\! } M \leq N_1 \cap N_2\mbox{, } \mathfrak{G}/M \iso F(X, W)\}.$$

The structure $\Gamma_{\mathfrak{G}, X} = (A^X_{\mathfrak{G}}, R^X_{\mathfrak{G}})$ is a graph, and moreover to any graph $\Gamma$ (and any Boolean space $X$) we may construct an Artin-Schreier structure $\mathfrak{G}_{\Gamma}$ such that $\Gamma_{\mathfrak{G}_{\Gamma},\, X} \iso \Gamma$; namely $\mathfrak{G}_{\Gamma} = F(X, G_{\Gamma})$ where $G_{\Gamma}$ is the group from \emph{Construction \ref{gconstruction}} ({this is \cite[p.\ 104, comment 1]{haran}}). Similarly, if $K$ is a field, define $A^{rc}_K$ to be the set of Galois extensions $L$ (contained in a fixed separable closure $K^s$ of $K$) such that $\mathfrak{Gal}(L/K) \iso F(X(K), D)$. Define the binary relation $R^{rc}_K$ on $A^{rc}_K$ by:
\begin{align*}
R^{rc}_K = \{(L_1, L_2) \in A^{rc}_K \times A^{rc}_K \mbox{ : } &L_1 \neq L_2 \mbox{ \&\ }
\exists N/K \mbox{ Galois s.t. } L_1 L_2 \subseteq N \\ 
&\mbox{\&\ } \mathfrak{Gal}(N/K) \iso F(X(K), W)\}.   
\end{align*}

This defines a graph $\Gamma^{rc}_K$. Finally, we define the \emph{canonical graph structure for $\mathfrak{G}$}, denoted $\Gamma_{\mathfrak{G}}$, as $\Gamma_{\mathfrak{G},\, X(\mathfrak{G})/G'}$. We have the following intuitive result:

\begin{lem}\label{221}
Let $G$ be a profinite group, $X$ a Boolean space, and $K$ a field. Then:
\begin{itemize}
    \item {There is a \emph{natural} isomorphism} $\Gamma_G \iso \Gamma_{F(X, G),\, X}$;
    \item $\Gamma^{rc}_K \iso \Gamma_{\mathfrak{G}_K,
    \,X(K)} = \Gamma_{\mathfrak{G}_K,\, X(\mathfrak{G}_K)/G_{K(\sqrt{-1})}} = \Gamma_{\mathfrak{G}_K}$;
    \item If $\phi: \mathfrak{H} \rightarrow \mathfrak{G}$ is a \emph{Frattini cover} of Artin-Schreier structures, then $\Gamma_{\mathfrak{G},\, X} \iso \Gamma_{\mathfrak{H},\, X}$.
\end{itemize}
\end{lem}

\proof
\emph{Comments 1, 2, \& 3} of \cite[p.\ 104]{haran}. \emph{Frattini covers of Artin-Schreier structures} are introduced and their basic properties proven in \cite[\S 1]{haran}. \bs

\subsection{Additional Tools}
Two tools needed for the main proof not already introduced are maximal pro-$\mathcal{C}$ quotients, and some notion of ``free product''. We will provide these ourselves here. 

Let $\mathcal{C}$ be a full formation of finite groups {containing the family of finite $2$-groups}, fix $\mathfrak{G} = \langle G, G', d: X(\mathfrak{G}) \rightarrow G\rangle$, define $\mathcal{N}$ to be the family of closed normal subgroups $N \triangleleft G$ such that $G/N \in \mathcal{C}$, and let $N_{\mathcal{C}} = \bigcap \mathcal{N}$. Then $G/N_{\mathcal{C}}$ is the maximal pro-$\mathcal{C}$ quotient of $G$. As $G' \leq G$ is open and of index $\leq 2$, $N_{\mathcal{C}} \leq G'$, hence we may form the (well-defined) quotient Artin-Schreier structure: 

$$\mathfrak{G}_{\mathcal{C}} \defeq \mathfrak{G}/N_{\mathcal{C}} = \langle G/N_{\mathcal{C}},\, G'/N_{\mathcal{C}},\, \overline{d} : X(\mathfrak{G})/N_{\mathcal{C}} \rightarrow G/N_{\mathcal{C}}\rangle.$$

This Artin-Schreier structure is defined to be the \emph{maximal pro-$\mathcal{C}$ quotient of $\mathfrak{G}$.} {Notice as $\mathcal{C}$ is a full formation, $G'/N_{\mathcal{C}}$ is pro-$\mathcal{C}$ (\cite[Lemma 17.3.1]{friedjarden}).}

\begin{lem}\label{lem1}
Suppose $\mathcal{C}$ is a full formation of finite groups {containing the family of finite 2-groups}, and $\mathfrak{G}$ is a projective Artin-Schreier structure. Then the maximal pro-$\mathcal{C}$ quotient $\mathfrak{G}_{\mathcal{C}}$ is also a projective Artin-Schreier structure.
\end{lem}

\proof
We will first show $\mathfrak{G}_{\mathcal{C}}$ is \textit{$\mathcal{C}$-projective} (projective relative to $\mathcal{C}$-groups; see \cite[Definition 22.3.1]{friedjarden}). Let $\mathfrak{A}, \mathfrak{B}$ be finite Artin-Schreier structures whose underlying groups $A, B$ are elements of $\mathcal{C}$ (immediately this implies $A', B' \in \mathcal{C}$). Consider the embedding problem
\[\begin{tikzcd}
& \mathfrak{G}_{\mathcal{C}} = \mathfrak{G}/N_{\mathcal{C}} \arrow[d, "\phi"]\\
\mathfrak{B} \arrow[r, "\alpha", twoheadrightarrow] & \mathfrak{A}
\end{tikzcd}\]

where $\alpha$ is an epimorphism and $\phi$ a morphism. Recall the cover that is the quotient morphism $\pi : \mathfrak{G} \twoheadrightarrow \mathfrak{G}/N_{\mathcal{C}}$ {(cf.\ \cite[pp.\ 458--459]{haranjarden})}; as $\mathfrak{G}$ is projective there exists a morphism $\gamma : \mathfrak{G} \rightarrow \mathfrak{B}$ such that $\alpha \circ \gamma = \phi \circ \pi$, {by \cite[Lemma 7.5]{haranjarden}}.

We claim $\gamma$ factors through $N_{\mathcal{C}}$. Indeed, $\gamma$ consists of the continuous morphisms $\gamma_1 : G \rightarrow B$ (which has the property $\gamma_1^{-1}(B') = G'$) and $\gamma_2 : X(\mathfrak{G}) \rightarrow X(\mathfrak{B})$ making the relevant diagrams commute. By the \textit{First Isomorphism Theorem for Groups}, $G/\ker(\gamma_1) \iso \gamma_1(G) \leq B$, hence as $\mathcal{C}$ is a full formation, by definition $\ker(\gamma_1) \in \mathcal{N}$. Therefore $N_{\mathcal{C}} \subseteq \ker(\gamma_1)$, and thus the map $\gamma_1$ factors through $N_{\mathcal{C}}$.

For $\gamma_2 : X(\mathfrak{G}) \rightarrow X(\mathfrak{B})$, notice that:
$$x_1 \sim x_2 \implies \gamma_2(x_2) = \gamma_2(x_1^{\sigma}) = \gamma_2(x_1)^{\gamma_1(\sigma)} = \gamma_2(x_1),$$

as $\sigma \in N_{\mathcal{C}} \subseteq \ker(\gamma_1)$. Thus, $\gamma_2$ factors through the equivalence relation $\sim$, and gives rise to a natural map $X(\mathfrak{G})/N_{\mathcal{C}} \rightarrow X(\mathfrak{B})$. Therefore the morphism $\gamma : \mathfrak{G} \rightarrow \mathfrak{B}$ factors through a morphism $\overline{\gamma} : \mathfrak{G}/N_{\mathcal{C}} \rightarrow \mathfrak{B}$, i.e.\ $\overline{\gamma} \circ \pi  = \gamma$. We conclude $\alpha \circ \overline{\gamma} = \phi$ as required to prove $\mathfrak{G}_{\mathcal{C}}$ is $\mathcal{C}$-projective. 

Now consider the below embedding problem where $\mathfrak{A}, \mathfrak{B}$ are arbitrary finite Artin-Schreier structures, and $\alpha$ is an epimorphism. We may assume WLOG $\phi$ is an epimorphism {(as Haran \& Jarden note in \cite[pp.\ 471--472]{haranjarden}, we may replace $\mathfrak{A}$ with an epimorphic image $\mathfrak{A}_0$ of $\mathfrak{G}_{\mathcal{C}}$, and $\mathfrak{B}$ with the fibred product of $\mathfrak{B}$ with $\mathfrak{A}_0$ over $\mathfrak{A}$. Cf.\ \cite[p.\ 460 \& Lemma 4.6]{haranjarden} for details on fibred products of Artin-Schreier structures)}.
\[
\begin{tikzcd}
& \mathfrak{G}_{\mathcal{C}} \arrow[d, "\phi"]\\
\mathfrak{B} \arrow[r, "\alpha", twoheadrightarrow] & \mathfrak{A}
\end{tikzcd}
\]

Hence $A, A' \in \mathcal{C}$. By\footnote{Ribes \& Zalesskii prove this for $\mathcal{C}$ a \emph{saturated variety of finite groups}; terminology found in \cite[pp.\ 19--20 \& Definition 7.6.4]{ribes}. We note a full formation of finite groups \emph{is} a `saturated variety', by definition (\cite[pp.\ 19--20]{ribes}) and \cite[Example 7.6.5 (1)]{ribes}.} \cite[Lemma 7.6.6]{ribes} there exists $M \leq B$ such that $B = \ker(\alpha) M$, $M \in \mathcal{C}$, and $\alpha(M) = A$. Let $M' = B' \cap M$. {We may further assume $\Inv(B) = \{\epsilon_1, \dots, \epsilon_n\} \subset M$, as if not the group $M\langle \epsilon_1 \rangle \dots \langle \epsilon_n \rangle$ has these properties (note we assumed $\mathcal{C}$ contains the family of finite $2$-groups).} Clearly $M' \leq M$ (and thus $M' \in \mathcal{C}$), but also $\alpha(M') = A'$. 

Indeed, $\alpha(M') \subseteq A'$, and furthermore given $a' \in A'$ we may find $b' \in B'$ such that $\alpha(b') = a'$. $M$ is constructed (\cite[Lemma 7.6.6]{ribes}) so that $\ker(\alpha) M = B$, hence $b' = k m$ for some $k \in \ker(\alpha)$ and $m \in M$. By definition of $\alpha$, $\ker(\alpha) \subseteq B'$, hence $m = k^{-1} b' \in B' \cap M$. Then $A' \subseteq \alpha(M')$ as desired. Note $M' \neq M \iff B' \neq B$, as $\ker(\alpha) M = B$ whereas $\ker(\alpha) M' \leq \ker(\alpha) B' = B'$. Finally, $(M : M') \leq (B : B') \leq 2$ as required to define the Artin-Schreier substructure
$$\mathfrak{M} = \langle M,\, M',\, X(\mathfrak{M}) \xrightarrow{d_M} M\rangle \quad \leq \quad \mathfrak{B} = \langle B,\, B',\,X(\mathfrak{B}) \xrightarrow{d_B} B\rangle,$$

where $X(\mathfrak{M}) \defeq d_B^{-1}(M)$ and $d_M = d_B|_{X(\mathfrak{M})}$. Note that as $B$ is finite, $M$ is closed, hence as $d_B$ is continuous the space $X(\mathfrak{M}) \subseteq X(\mathfrak{B})$ is Boolean. Therefore $\mathfrak{M}$ indeed satisfies the definition of an Artin-Schreier structure. It is furthermore a \emph{substructure} of $\mathfrak{B}$, as the canonical inclusion map $\mathfrak{M} \hookrightarrow \mathfrak{B}$ is a morphism, {and by construction $\alpha : \mathfrak{M} \twoheadrightarrow \mathfrak{A}$ is an epimorphism.}

As $\mathfrak{G}_{\mathcal{C}}$ is $\mathcal{C}$-projective, there exists a morphism $\gamma : \mathfrak{G}_{\mathcal{C}} \rightarrow \mathfrak{M}$ which solves the embedding problem. Solving finite embedding problems is sufficient by \cite[Lemma 7.5]{haranjarden} to ensure $\mathfrak{G}_{\mathcal{C}}$ is projective in the category of Artin-Schreier structures, as required. \bs

Examining \emph{Theorem \ref{pacfiniteaxiom}} reveals in fact we do not need a \textit{coproduct} in the category of Artin-Schreier structures --  we shall see that the following construction will suffice (cf.\ \cite[\S 3.3]{fehmthesis} for a similar technique).

Let $\mathfrak{G} = \langle G,\, G',\, d: X(\mathfrak{G}) \rightarrow G\rangle$ be an Artin-Schreier structure and $E$ a profinite group. Define $\mathfrak{G} \Diamond E := \langle G \star E,\, G'_{\star},\, \Inv(G \star E) \hookrightarrow G \star E\rangle$, where ``$\star$'' is the coproduct in the category of profinite groups, $G'_{\star} := \ker(G \star E \rightarrow \Z\!/2\Z)$ where $G \star E \rightarrow \Z\!/2\Z$ is uniquely determined by $E \rightarrow \{1\}$, $G \rightarrow G/G'$, {where $G/G' \iso \Z\!/2\Z$ or $\{1\}$}, and $\Inv(G \star E)$ is the set of involutions of $G \star E$. Note that a priori this is not necessarily even a weak Artin-Schreier structure. Though a posteriori, we have the following result.

\begin{lem}\label{lem2}
Let $\mathfrak{G}$ be a projective Artin-Schreier structure, and $E$ a projective profinite group. Then $\mathfrak{G} \Diamond E$ is a projective Artin-Schreier structure.
\end{lem}

\proof
(Cf.\ \emph{Lemma \ref{lem3}} for a different style of argument.) By \cite[Theorem 2.1]{harancohom2}, $G'$ is a projective profinite group and $X(\mathfrak{G}) \iso \mbox{Inv}(G)$. Note that $G'_{\star}$ is open, normal, and of index $\leq 2$ in $G \star E$. Therefore in order for $\mathfrak{G} \Diamond E$ to be a weak Artin-Schreier structure, it is simply required that $\Inv(G \star E)$ be a closed subset of $G \star E$ (as then it is Boolean). Equivalently, by \cite[Remark 7.6]{haranjarden} we wish there to exist an open $U \triangleleft G \star E$ such that $U \cap \Inv(G \star E) = \emptyset$. This is satisfied by $U = G'_{\star}$; {indeed, if $\epsilon \in G \star E$ is an involution, by \cite[Theorem A]{herfort} it is conjugate to an involution $\sigma$ of $G$ or $E$. However $G' \leq G$ and $E$ are projective, hence torsion free, so $\sigma \in G \setminus G'$. Therefore $\epsilon \not\in \ker(G \star E \rightarrow \Z/2\Z)$ above, thus $\epsilon \not\in G'_{\star}$.}

We shall show every finite embedding problem for $\mathfrak{G} \Diamond E$ has a solution, which completes the proof by \cite[Lemma 7.5]{haranjarden}. Consider the following embedding problem:
\[
\begin{tikzcd}
& \mathfrak{G} \Diamond E \arrow[d, "\varphi"]\\
\mathfrak{B} \arrow[r, "\alpha", twoheadrightarrow] & \mathfrak{A}
\end{tikzcd}
\]

where $\varphi$ is a morphism and $\alpha$ is an epimorphism of finite weak Artin-Schreier structures $\mathfrak{A}, \mathfrak{B}$. By \cite[Lemma 7.5]{haranjarden}, as the forgetful map of $\mathfrak{G}\Diamond E$ is injective it will suffice to assume $X(\mathfrak{A}) \subseteq A$, $X(\mathfrak{B}) \subseteq B$. Diagrammatically:

\[
\begin{tikzcd}[column sep = tiny, row sep = tiny]
& \mathfrak{G} \arrow[ddddd, "\varphi", start anchor={[xshift=-2.5mm]}, end anchor={[xshift=-2.5mm]}]\Diamond E \hspace{2mm}\arrow[phantom, r, "=", pos=0.37, scale=1] & \langle G \star E \arrow[ddddd, "\varphi", swap], & G'_{\star}, \arrow[ddddd, "\varphi", swap]&  \Inv(G\star E) \arrow[ddddd, "\varphi", swap] \arrow[rrr, "incl.", hook] \arrow[rrr, ""{name=U, below}]{} & & & G \star E \arrow[ddddd, "\varphi"] \rangle\\
&&&&&&&\\
&&&&&&&\\
&&&&&&&\\
&&&&&&&\\
\mathfrak{B} \arrow[ddddddd, phantom, sloped, anchor=center, "=", scale=1.2] \arrow[r, "\alpha", twoheadrightarrow] &  \mathfrak{A}\hspace{2mm} \arrow[phantom, r, "=", pos=0.43, scale=1] &\langle A, & A',&X(\mathfrak{A}) \arrow[rrr, phantom, ""{name=D, above}]{} \arrow[rrr,phantom, ""{name=U2, below}]{} \arrow[rrr, "incl.", hook] \arrow[phantom, from=U, to=D, "\circlearrowleft" scale=1.5] & &
&  A\rangle\\
&&&&&&&\\
&&&&&&&\\
&&&&&&&\\
&&&&&&&\\
&&&&&&&\\
&&&&&&&\\
\langle B \arrow[rruuuuuuu, start anchor={[xshift=-1mm]}, "\alpha", twoheadrightarrow],  & \hspace{6mm}B'  \arrow[rruuuuuuu, "\alpha", pos=0.48, twoheadrightarrow, start anchor={[xshift=-3.5mm, yshift=-1mm]}],\hspace{7mm} &  X(\mathfrak{B}) \arrow[rruuuuuuu, "\alpha", pos=0.49, twoheadrightarrow, start anchor={[xshift=-2mm]}] \arrow[rr, phantom, ""{name=D2, above}]{} \arrow[rr, "incl.", hook] & & B \arrow[rrruuuuuuu, "\alpha", twoheadrightarrow, pos=0.66, swap, start anchor={[xshift=-3mm]}] \arrow[phantom, from=U2, to=D2, "\circlearrowleft" scale=1] \rangle &&&\\
\end{tikzcd}
\]

(Note the maps are indeed correctly labelled; each map is $\varphi: G \star E \rightarrow A$, or $\alpha : B \rightarrow A$, restricted to a subset of the relevant group.) This diagram gives rise to:
\[
\begin{tikzcd}
& G \arrow[d, "\varphi|_G"] &&&&E \arrow[d, "\varphi|_E"]\\
B \arrow[r, "\alpha", twoheadrightarrow] & A&&&B \arrow[r, "\alpha", twoheadrightarrow] & A
\end{tikzcd}
\]

where $\alpha : B' \twoheadrightarrow A
'$. The former is solvable by $\gamma_G : G \rightarrow B$, as $G$ is real projective {(\cite[Proposition 7.7]{haranjarden})} and this is a \emph{real} embedding problem for $G$. {Indeed (following the exposition of \cite[p.\ 474]{haranjarden}) we may assume in addition that $X(\mathfrak{B}) = \Inv(B \setminus B')$, as by \cite[Corollary 6.2]{haranjarden} there exists a finite group $B_1$ and epimorphism $\theta : B_1 \twoheadrightarrow B$ such that the image under $\theta$ of $\Inv(B_1 \setminus \ker(\theta))$ is exactly $X(\mathfrak{B})$, and we may replace $\mathfrak{B}$ with $\langle B_1,\, \theta^{-1}(B'),\, \Inv(B_1 \setminus \ker(\theta)) \hookrightarrow B_1 \rangle$ if desired. Therefore if $\epsilon \in G$ is an involution with $\varphi|_G(\epsilon) \neq 1$, considering $\epsilon$ as an involution of $G\star E$, under $\alpha : X(\mathfrak{B}) \twoheadrightarrow X(\mathfrak{A})$ there exists an involution $b \in B \setminus B'$ with $\alpha(b) = \varphi(\epsilon) = \varphi|_G(\epsilon)$, as required.}

The latter diagram is solvable by $\gamma_E : E \rightarrow B$ (as $E$ is projective). By definition of the free product, there exists a morphism $\gamma : G \star E \rightarrow B$ {uniquely extending $\gamma_G$, $\gamma_E$} such that $\alpha \circ \gamma = \varphi$. 

Consider $\gamma|_{G'_{\star}} : G'_{\star} \rightarrow B$. We wish that $\gamma|_{G'_{\star}}(G'_{\star}) \subseteq B'$. This is indeed the case: by {definition} $\alpha(B') = A'$ (and moreover $\alpha(B \setminus B') \subseteq A \setminus A'$), so since $\alpha \circ \gamma = \varphi : G'_{\star} \rightarrow A'$ it must be that $\gamma|_{G'_{\star}}(G'_{\star}) \subseteq B'$ as desired. We may update the main diagram as follows:

\[
\begin{tikzcd}[column sep = tiny, row sep = tiny]
& &\mathfrak{G} \arrow[ddddd, "\varphi", start anchor={[xshift=-2.5mm]}, end anchor={[xshift=-2.5mm]}, pos=0.4]\Diamond E \hspace{5mm}\arrow[phantom, rr, "=", pos=0.37, scale=1] & & \langle G \star E \arrow[ddddd, "\varphi", swap, pos=0.4] \arrow[dddddddddddddllll, end anchor={[xshift=-1mm]}, bend right=9, pos=0.4, swap, red, "\gamma"],& & G'_{\star}, \arrow[ddddd, "\varphi", pos=0.4]& & \Inv(G\star E) \arrow[ddddd, "\varphi", swap, pos=0.4] \arrow[rrr, "incl.", hook] \arrow[rrr, ""{name=U, below}]{} & & & G \star E \arrow[ddddd, "\varphi",pos=0.4] \rangle\\
&&&&&&&&&&&\\
&&&&&&&&&&&\\
&&&&&&&&&&&\\
&&&&&&&&&&&\\
\mathfrak{B} \arrow[dddddddd, phantom, sloped, anchor=center, "=", scale=1.2] \arrow[rr, "\alpha", twoheadrightarrow] & &  \mathfrak{A}\hspace{5mm} \arrow[phantom, rr, "=", pos=0.4, scale=1] & &\langle A,& & A',& &X(\mathfrak{A}) \arrow[rrr,phantom, ""{name=D, above}]{} \arrow[rrr,phantom, ""{name=U2, below}]{} \arrow[rrr, "incl.", hook] \arrow[phantom, from=U, to=D, "\circlearrowleft" scale=1.5] & &
&  A\rangle\\
&&&&&&&&&&&\\
&&&&&&&&&&&\\
&&&&&&&&&&&\\
&&&&&&&&&&&\\
&&&&&&&&&&&\\
&&&&&&&&&&&\\
&&&&&&&&&&&\\
\langle B \arrow[rrrruuuuuuuu, start anchor={[xshift=0mm]}, "\alpha", twoheadrightarrow], & & \hspace{6mm}B' \arrow[from=uuuuuuuuuuuuurrrr, swap, crossing over, bend right=9, end anchor={[xshift=-3mm]}, red, pos=0.4, "\gamma"] \arrow[rrrruuuuuuuu, "\alpha", pos=0.48, twoheadrightarrow, start anchor={[xshift=-6mm, yshift=-1mm]}],\hspace{7mm} &  X(\mathfrak{B}) \arrow[rrrrruuuuuuuu, "\alpha", pos=0.51, twoheadrightarrow] \arrow[r, phantom, ""{name=D2, above}]{} \arrow[r, "incl.", swap, hook, end anchor={[xshift=9mm]}] &   \hspace{10mm}B \arrow[rrrrrrruuuuuuuu, "\alpha",pos=0.65, twoheadrightarrow, swap] \arrow[phantom, from=U2, to=D2, "\circlearrowleft" scale=1] \arrow[from=uuuuuuuuuuuuurrrrrrr, swap, bend right=12, end anchor={[xshift=1.5mm]}, crossing over, red, "\gamma"]\rangle &&&&&&&\\
\end{tikzcd}
\]

If it is the case that $\gamma(\Inv(G \star E)) \subseteq X(\mathfrak{B})$, then we are finished. {Recall $X(\mathfrak{B}) = \Inv(B \setminus B')$; then} for all $\varepsilon \in \Inv(G \star E)$, since $\alpha(\gamma(\varepsilon)) = \varphi(\varepsilon) \in X(\mathfrak{A})$, by the initial set up $\gamma(\varepsilon) \in B \setminus B'$ and hence $\gamma(\varepsilon) \in X(\mathfrak{B})$ as desired.

We conclude that $\gamma : \mathfrak{G} \Diamond E \rightarrow \mathfrak{B}$ is a morphism of weak Artin-Schreier structures solving the initial finite real embedding problem. This completes the proof. \bs

\noindent We will reference the following from the proof of \cite[Theorem 3.1]{haran}:

\begin{lem}\label{arnomap}
Let $K$ be a PRC field; $\Gamma^{rc}_K$ is interpretable $($in the sense of \emph{Definition \ref{interp}}$)$ in $K$.
\end{lem}

\proof
This is \cite[pp.\ 106--107]{haran}, which we elaborate on now. Let $G$ be a finite group, and $\alpha_{r, G}(x_1, \dots x_r)$ be an $\mathcal{L}_r$-formula such that
\begin{align*}
K \models \alpha_{r, G}(\overline{a}) \iff\,&\mbox{$f_{\overline{a}}(T) = T^r + a_1 T^{r-1} + \dots + a_r$ is irreducible over $K$, and}\\
&\mathfrak{Gal}(K_{\overline{a}}(\sqrt{-1})/K) \iso F(X(K), G),    
\end{align*}

where $\overline{a} \in K^r$. Indeed, $\alpha_{r, G}(\overline{a})$ can be given by the conjunction of the following statements:
\begin{itemize}
\item $f_{\overline{a}}$ is irreducible over $K$, $K_{\overline{a}}/K$ is Galois, and $\Gal(K_{\overline{a}}/K) \iso G$;
\item $\sqrt{-1} \not\in K$;
\item $K_{\overline{a}}/K$ is {totally real}.\label{totreal}
\end{itemize}

All but the last are standard to express in the language of rings, and the last is covered by Prestel in \cite[Theorem 4.1]{prestel}. On \emph{p.\ 154 ibid.}, Prestel shows an equivalent formulation of ``$K_{\overline{a}}/K$ is totally real'' is $\mathcal{L}_r$-axiomatisable, assuming $K$ is PRC. He assumes (in our notation) that $f_{\overline{a}}$ is absolutely irreducible, though this is not used in this part of his proof. 

Fix finite groups $D, W$ satisfying the graph conditions (\emph{Construction \ref{gconstruction}}), and define $l = |D|$. As in \emph{Construction \ref{anothcons}} one may construct an $\mathcal{L}_r$-formula $\rho_{D, W}$ such that for $\overline{b}, \overline{c} \in K^l$,}
\begin{align*}
K \models \rho_{D, W}(\overline{b}, \overline{c}) \iff \,& K_{\overline{b}} \not\subseteq K_{\overline{c}} \lor K_{\overline{c}} \not\subseteq K_{\overline{b}}\mbox{, \hspace{2mm}} K \models \alpha_{l, D}(\overline{b}) \land \alpha_{l, D}(\overline{c})\mbox{, \hspace{2mm}and}\\
&{K \models \exists \overline{z}\,( \alpha_{|W|, W}(\overline{z})\land \mbox{``}K_{\overline{b}} \subseteq K_{\overline{z}}\mbox{''} \land \mbox{``}K_{\overline{c}} \subseteq K_{\overline{z}}\mbox{''}).}    
\end{align*}

We may then define a recursive translation map $-' : \mbox{Form}(\mathcal{L}_{gr}) \rightarrow \mbox{Form}(\mathcal{L}_r)$; $\phi \mapsto \phi'$ by the following rules:
\begin{itemize}
    \item $R(X, Y) \mapsto (R(X, Y))' = \rho_{D, W}(\overline{x}, \overline{y})$;
    \item $\neg \varphi \mapsto \neg(\varphi')$; 
    \item $\varphi_1 \land \varphi_2 \mapsto (\varphi_1') \land (\varphi_2')$;
    \item $\exists x (\varphi) \mapsto \exists \overline{x}\, (\alpha_{l, D}(\overline{x}) \land \varphi'(\overline{x}))$. 
\end{itemize}

Recalling \emph{Definition \ref{interp}}, set $\mathcal{L}_1 = \mathcal{L}_r$, $\mathcal{L}_0 = \mathcal{L}_{gr}$, $n = l= |D|$, $\delta(\overline{x}) = \alpha_{n, D}(\overline{x})$, and $f : \alpha_{n, D}(K^n) \rightarrow \Gamma^{rc}_{K}$; $\overline{a} \mapsto K_{\overline{a}}(\sqrt{-1})$.  This is indeed surjective: by \emph{Remark \ref{snaky}} if $L \in A_K^{rc}$ there exists a totally real Galois extension $L_0/K$  such that $L = L_0(\sqrt{-1})$ and $\mbox{Gal}(L_0/K) \iso D$, hence $L_0$ is the splitting field of a degree $n = |D|$ monic separable irreducible polynomial $f_{\overline{a}}(T) = T^n + a_1 T^{n-1} + \dots + a_n$ over $K$. Therefore $K \models \alpha_{n, D}(\overline{a})$.

Note \emph{Definition \ref{interp} (2)} is satisfied by the above construction of $-'$, and condition $(\dagger)$ is confirmed in \cite[Theorem 3.1]{haran}. \bs

\begin{remark}
Notice the $\mathcal{L}_r$-formula $\alpha_{n,D}(\overline{x})$ and the map $-' : \mbox{Form}(\mathcal{L}_{gr}) \rightarrow \mbox{Form}(\mathcal{L}_r)$ of \emph{Lemma \ref{arnomap}} do not depend on $K$ or $\Gamma^{rc}_K$; the interpretation is \emph{uniform} across PRC fields $K$.\sq
\end{remark}

\subsection{Results} Fix a PRC field $K$ with a finite subtheory $\Sigma \subseteq \Th(K;\mathcal{L}_r)$. {By the \emph{Compactness Theorem}, there exists a finite set of $\mathcal{L}_r$-sentences $\Delta$ such that $\Delta \models \Sigma$ and $\Delta = \Delta_1 \cup \Delta_2 \cup \Delta_3 \cup \Delta^*_4$, where $\Delta_1$ is a finite subset of $\mathcal{L}_r$-sentences specifying the characteristic and degree of imperfection of $K$ (\emph{Corollary \ref{prcax} (1)}), $\Delta_2$ is a finite subset of {\ttfamily PRC} (\emph{Corollary \ref{prcax} (2)}), $\Delta_3$ is a finite subset of $\Th^{alg}(K)$ (\emph{Corollary \ref{prcax} (3)}), $\Delta_4^*$ is a finite subset of $\Th(S(G_K); \mathcal{L}_G(S(G_{K_0})))^*$ (\emph{Corollary \ref{prcax} (4)}), and $\Delta_4$ is a finite subset of $\Th(S(G_K); \mathcal{L}_G(S(G_{K_0})))$ with $\varphi \in \Delta_4 \Leftrightarrow \varphi^* \in \Delta_4^*$.

Let $\Lambda$ be the set of universal sentences of $\Delta_3$, and let $\Kbad$ be the join of $K(\sqrt{-1})$ and of minimal Galois extensions $F/K$ within a fixed algebraic closure $\widetilde{K}$ of $K$, with $F \models \neg \lambda$ for $\lambda \in \Lambda$. Let $\overline{a}_{\Sigma} \in S(G_{K_0}) \subset S(G_K)$ be a finite tuple of elements such that $\Delta_4$ is a set of finitely many $\mathcal{L}_G(\overline{a}_{\Sigma})$-sentences. Fix $n_{\Sigma} \in \N$ such that $S_1, \dots, S_{n_{\Sigma}}$ is the smallest consecutive sequence of sorts involving the sentences of $\Delta_4$. Let $\widehat{p}$ be the smallest odd prime larger than $n_{\Sigma} + |\Gal(\Kbad/K)|$, $P_{\Sigma}$ the set of primes $\{2, 3, \dots, \widehat{p}\}$, and $\mathcal{C}_{\Sigma}$ the formation of finite groups whose order is necessarily a product of powers of primes of $P_{\Sigma}$ (including trivial powers).} We have the following theorem:

\begin{thm}\label{mainthm2}
Assume the above setup. For any nontrivial graph $\Gamma$, there exists a PRC field $K_{\Gamma} \supseteq K$ such that 
\begin{enumerate}
    \item $K_{\Gamma} \models \Sigma$;
    \item $\Gamma^{rc}_{K_{\Gamma}} \iso \Gamma$;
    \item $X(K_{\Gamma}) \iso X(K)$;
\end{enumerate}
\end{thm}

\proof
Assume $K$ is formally real; if it is not the result is \emph{Theorem \ref{thisolthm}}. We will mirror \emph{Theorem \ref{thisolthm}--Lemma \ref{graphimp}} closely to prove (1) \&\ (2).

{Recall from \emph{Construction \ref{graph}} the groups $D_{\Sigma}$, $W_{\Sigma}$. Fix a nontrivial graph $\Gamma$, and consider the profinite group $G_{\Gamma}$ from \emph{Construction \ref{gconstruction}}. Let $\widetilde{G_{\Gamma}}$ be the universal Frattini cover of $G_{\Gamma}$ with map $f : \widetilde{G_{\Gamma}} \rightarrow G_{\Gamma}$. By \cite[Corollary 54]{cvddmac} the sorts $S_1, \dots, S_{\widehat{p}}$ of $S(\widetilde{G}_{\Gamma})$ are trivial (there is no proper normal subgroup $N \triangleleft \widetilde{G_{\Gamma}}$ with $[\widetilde{G_{\Gamma}} : N] \leq \widehat{p}$), and $\Gamma \iso \Gamma_{\widetilde{G_{\Gamma}}}$ by \cite[Lemma 28.6.1]{friedjarden}.} Consider further the diamond product $\mathfrak{G}^* = \mathfrak{G}_{P_{\Sigma}} \Diamond \widetilde{G_{\Gamma}}$, where $\mathfrak{G}_{P_{\Sigma}} = \langle G_{P_{\Sigma}}, G'_{P_{\Sigma}}, d: X(\mathfrak{G}_{P_{\Sigma}}) \rightarrow G_{P_{\Sigma}}\rangle$ is the maximal pro-$\mathcal{C}_{\Sigma}$ quotient of $\mathfrak{G}_K$ (\emph{Lemma \ref{lem1}}). $\mathfrak{G}^*$ is projective by \emph{Lemma \ref{lem2}}, {and $S(G_{P_{\Sigma}}),\, S(G_{P_{\Sigma}} \star \widetilde{G_{\Gamma}}) \models \Delta_4$ by \emph{Lemma \ref{addedsea}}.}

We claim $\Gamma_{G'_{\star}} \iso \Gamma$. We (almost) repeat \emph{Lemma \ref{graphimp}:} note $\widetilde{G_{\Gamma}}, G'_{P_{\Sigma}} \leq G'_{\star}$, and 
$$\widetilde{G_{\Gamma}} \iso (G_{P_{\Sigma}} \star \widetilde{G_{\Gamma}}) / \langle G_{P_{\Sigma}} \rangle \iso G'_{\star}/(G'_{\star} \cap \langle G_{P_{\Sigma}} \rangle).$$

(The second isomorphism results from the \textit{Second Isomorphism Theorem for Groups}.) Hence there is an epimorphism $\pi : G'_{\star} \rightarrow \widetilde{G_{\Gamma}}$ which yields an embedding $\Gamma_{\widetilde{G_{\Gamma}}} \hookrightarrow \Gamma_{G'_{\star}}$. Conversely, if $M \triangleleft G'_{\star}$ is open and $G'_{\star}/M \iso D_{\Sigma}$ (resp.\ $W_{\Sigma}$), under $\rho : G'_{\star} \twoheadrightarrow D_{\Sigma}$ (resp.\ $W_{\Sigma}$) we have $\rho(G'_{P_{\Sigma}}) \leq D_{\Sigma}$ (resp.\ $W_{\Sigma}$). As $G'_{P_{\Sigma}}$ is pro-$\mathcal{C}_{\Sigma}$ (\cite[Lemma 17.3.1]{friedjarden}), any quotient of $G'_{P_{\Sigma}}$ is a $P_{\Sigma}$-group, hence there exists elements of $\rho(G'_{P_{\Sigma}})$ of the wrong order unless $G'_{P_{\Sigma}} \leq \ker(\rho)$. Hence $\ker(\pi) \leq \ker(\rho)$, thus $M = \pi^{-1}(\pi(M))$, and therefore $\pi^{-1}$ induces an isomorphism of graphs $\Gamma_{G'_{\star}} \iso \Gamma$ as claimed.

The canonical graph structure for $\mathfrak{G}^*$ is also recovered correctly: 
\begin{align*}
  \Gamma_{\mathfrak{G}^*} &= \Gamma_{\mathfrak{G}^*,\, X(\mathfrak{G}^*)/G'_{\star}} &&\mbox{by definition,}\\
  &\iso \Gamma_{F(X(\mathfrak{G}^*)/G'_{\star},\, G'_{\star}),\, X(\mathfrak{G}^*)/G'_{\star}}&&\mbox{by definition,}\\  
  &\iso \Gamma_{F(X(\mathfrak{G}_{P_{\Sigma}})/G'_{P_{\Sigma}},\, G'_{\star}),\, X(\mathfrak{G}_{P_{\Sigma}})/G'_{P_{\Sigma}}}&&\mbox{by \emph{Lemma \ref{221}},}\\
  &\iso \Gamma_{G'_{\star}}&&\mbox{by \emph{Lemma \ref{221},}} \\
  &\iso \Gamma &&\mbox{by the above}.
\end{align*}

\label{arghere}There is an epimorphism of Artin-Schreier structures $\mathfrak{G}_{P_{\Sigma}} \twoheadrightarrow \mathfrak{Gal}(\Kbad/K)$. Indeed, as $G_{P_{\Sigma}} = G_K/N_{\mathcal{C}_{\Sigma}}$, if $E_{\mathcal{C}_{\Sigma}}$ is the fixed field $(\widetilde{K})^{N_{\mathcal{C}_{\Sigma}}}$ then we have the tower of fields $\widetilde{K}/E_{\mathcal{C}_{\Sigma}}/\Kbad/K$ by design. \cite[Example 3.4 (a)]{haranjarden} ensures there is an epimorphism $\mathfrak{G}_{P_{\Sigma}} \twoheadrightarrow \mathfrak{Gal}(\Kbad/K)$. By \emph{Corollary \ref{214}} there is an extension $K_{\Gamma}$ of $K$ that is PRC, $K_{\Gamma} \cap \Kbad = K$, and $\mathfrak{G}_{K_{\Gamma}} \iso \mathfrak{G}^*$. Thus, $S(G_{K_{\Gamma}}) \models \Delta_4$, and $K_{\Gamma} \models \Delta_3$. Therefore $K_{\Gamma} \models \Delta_1 \cup \Delta_2 \cup \Delta_3 \cup \Delta^*_4$ (implying $K_{\Gamma} \models \Sigma$), {and from \emph{Lemma \ref{221}},}
$$\Gamma^{rc}_{K_{\Gamma}} \iso \Gamma_{\mathfrak{G}_{K_{\Gamma}},\, X(K_{\Gamma})} \iso \Gamma_{\mathfrak{G}^*,\, X(\mathfrak{G}^*)/G'_{\star}} \iso \Gamma,$$

This proves (1) \&\ (2). The proof of (3) is a direct consequence of the above setup: {$X(K_{\Gamma}) \iso X(\widetilde{K_{\Gamma}}/K_{\Gamma})/G_{K_{\Gamma}(\sqrt{-1})} \iso X(\mathfrak{G}^*)/G'_{\star}$, as $\mathfrak{G}_{K_{\Gamma}} \iso \mathfrak{G}^*$. Furthermore:}
\begin{align*}
    X(\mathfrak{G}^*)/G'_{\star}&= \Inv(G_{P_{\Sigma}} \star \widetilde{G_{\Gamma}})/G'_{\star}\\ 
    &\iso \Inv(G_{P_{\Sigma}})/G'_{P_{\Sigma}} &&\mbox{from\footnotemark{} \cite[Theorem A]{herfort},}\\ 
    &\iso X(\mathfrak{G}_{P_{\Sigma}})/G'_{P_{\Sigma}} &&\mbox{as $\mathfrak{G}_{P_{\Sigma}}$ is projective; \cite[Proposition 7.4]{haranjarden},}\\
    &= X(\mathfrak{G}_K / N_{\mathcal{C}_{\Sigma}})/(G'_K / N_{\mathcal{C}_{\Sigma}})\\
    &= (X(\mathfrak{G}_K)/N_{\mathcal{C}_{\Sigma}})/(G'_K / N_{\mathcal{C}_{\Sigma}}).
\end{align*}

\footnotetext{The involutions of $G_{P_{\Sigma}} \star \widetilde{G_{\Gamma}}$ are exactly the conjugates of $\Inv(G_{P_{\Sigma}})$ in $G_{P_{\Sigma}} \star \widetilde{G_{\Gamma}}$, by \cite[Theorem A]{herfort}. Therefore the morphism $\mathfrak{G}^* \rightarrow \mathfrak{G}_{P_{\Sigma}}$ induces a continuous \emph{bijection} $\Inv(G_{P_{\Sigma}} \star \widetilde{G_{\Gamma}})/G'_{\star} \rightarrow \Inv(G_{P_{\Sigma}})/G'_{P_{\Sigma}}$.}Notice $X(\mathfrak{G}_{K})/G'_{K} \iso (X(\mathfrak{G}_K)/N_{\mathcal{C}_{\Sigma}})/(G'_K / N_{\mathcal{C}_{\Sigma}})$ as $\mathfrak{G}_K \twoheadrightarrow \mathfrak{G}_K/N_{\mathcal{C}_{\Sigma}}$ is a cover. Finally, $X(\mathfrak{G}_{K})/G'_{K} \iso X(K)$ ensuring $X(K_{\Gamma}) \iso X(K)$; the order space of $K_{\Gamma}$ is up to homeomorphism that of $K$. \bs

\begin{cor}\label{mainprcfieldundeccor}
Every PRC field is finitely undecidable.
\end{cor}\enlargethispage*{0.5\baselineskip}

\proof
Fix $K$ a PRC field and $\Sigma \subseteq \Th(K; \mathcal{L}_r)$ a finite subtheory. Assume $K$ is formally real (otherwise the result is \emph{Corollary \ref{mainpacfieldundeccor}}). Let $\Gamma$ be a nonempty graph; by \emph{Lemma \ref{arnomap} \& Theorem \ref{mainthm2}} there exists a PRC field $K_{\Gamma} \models \Sigma$ and $\Gamma$ is interpretable in $K_{\Gamma}$. Furthermore, the interpretation of \emph{Lemma \ref{arnomap}} is uniform in the sense of \emph{Definition \ref{uniformparam}}; the class of nontrivial graphs is uniformly interpretable in the class of PRC fields satisfying $\Sigma$. We conclude {\ttfamily PRC}$\,\cup\, \Sigma$ is hereditarily undecidable by \emph{Corollary \ref{endcor}}, hence $\Sigma$ is undecidable as required. \bs

\begin{cor}\label{prcnotfinax}
No PRC field is finitely axiomatisable. \bs
\end{cor}

\begin{remark}\label{theopenqs}
In 2014, Shlapentokh \& Videla \cite[\S 6]{shlapentokhvidela} posed three open questions about finite undecidability (recall their term is `\textit{finite hereditary} undecidability'):
\begin{enumerate}
    \item \textit{``It is known that the theory of the field of all totally real algebraic numbers is decidable \dots Is this theory finitely hereditarily undecidable?''}
    \item \textit{``Is the theory of pseudo real closed fields \dots finitely hereditarily undecidable?''}
    \item \textit{``In general, if $\mathfrak{T}$ is a theory of any subfield of $\widetilde{\Q}$, the algebraic closure of $\Q$, is $\mathfrak{T}$ finitely hereditarily undecidable?''}
\end{enumerate}
{(Quotes from \cite[p.\ 1262]{shlapentokhvidela}.)} We may answer these as follows:
\begin{enumerate}
    \item By Pop \cite{popprc}, the field $\Q^{tr}$ of totally real algebraic numbers is PRC, hence by \emph{Corollary \ref{mainprcfieldundeccor}} the answer to (1) is \emph{yes}.
    \item As posed, this is covered by Haran \cite{haran}, however \emph{Corollary \ref{mainprcfieldundeccor}} answers {in the positive} the more general subsequent question: \textit{is the theory of any pseudo real closed field finitely undecidable?} (Equivalently: \textit{is any completion of the theory of PRC fields finitely undecidable?})
    \item We cannot answer this in full generality, though we make the following comment: as $\Gal(\widetilde{\Q}/\!\Q)$ is profinite, it has a unique Haar measure $\mu$ \cite[Prop. 18.2.1]{friedjarden}. By the PAC Nullstellensatz \cite[Theorem 2.5]{jardenpac}, the fixed field $\widetilde{\Q}(\sigma_1, \dots, \sigma_n)$ is PAC for a $\mu$-measure 1 subset of $\Gal(\widetilde{\Q}/\!\Q)^n$. By \emph{Corollary \ref{mainpacfieldundeccor}}, we conclude finite undecidability `for almost all' fixed fields $\widetilde{\Q}(\sigma_1, \dots, \sigma_n)$. Of course, with \emph{Corollary \ref{mainprcfieldundeccor}}, we know there are even more infinite algebraic extensions $R/\!\Q$ with $R$ finitely undecidable, namely the formally real PRC fields $R$. \sq
\end{enumerate}
\end{remark}

\begin{remark}
The move from PAC to PRC is also quite natural from the perspective of classification theory. Recall the notion of a \textit{(super)rosy theory} (see \cite[Fact 4.4]{rosy}, and \cite{krupinski, onshuus} for further background). The following is \cite[Conjecture 3]{krupinski}:

\begin{customconj}{(Superrosy Fields)}
\textit{Every infinite superrosy field is perfect, bounded, and PRC.}\bsn
\end{customconj}

\begin{cor}
Assume the Superrosy Fields Conjecture. Then every infinite superrosy field is finitely undecidable. \bsn
\end{cor}

Also note that we can recover the specific proof of \emph{Theorem \ref{pacfiniteaxiom}} in the PRC field context -- that is, we do not need to conclude this using finite undecidability. The ingredients of this proof (and that of \emph{Theorem \ref{pacfiniteaxiom}}) are considerably more accessible than the ingredients of \emph{Corollary \ref{mainprcfieldundeccor}} (resp.\ \emph{Corollary \ref{mainpacfieldundeccor}}), hence we may speculate this argument could work in a considerably broader context. For brevity we do not include this here; this will be generalised in future work. \sqn
\end{remark}

\bigskip
\section{Pseudo-\texorpdfstring{$p$}{}-Adically Closed Fields}\label{stuffforppc}

The definition of a \emph{pseudo-$p$-adically closed field} is the natural step after that of PAC and PRC fields, from both an algebraic and model-theoretic standpoint. We shall assume the reader is familiar with these fields, and the paper \cite{haranjardenppc}.

One might naturally expect the case of P$p$C fields to be more intricate than what followed previously; the approach of this thesis is through the absolute Galois group of such fields, and even for $\Q_p$, $\Th(S(G_{\Q_p});\mathcal{L}_G)$ is monstrous compared to $\Th(S(G_{K});\mathcal{L}_G)$ for $K$ algebraically, real, or separably closed. A paper of Efrat recovers Haran's result:

\begin{thm}
\emph{\cite[Corollary 4.3]{efrat}.} The theory of formally $p$-adic P$p$C fields is (hereditarily) undecidable. \bsn
\end{thm}

We are unable to prove the finite undecidability of any given P$p$C field (see \emph{Remark \ref{lackofppc}}). In line with \emph{Theorem \ref{pacfiniteaxiom} \& Corollary \ref{prcnotfinax}}, we are able to show:

\begin{customthm}{\ref{noppcfinax}}
\textit{No bounded P$p$C field is finitely axiomatisable (even among the class of P$p$C fields).}
\end{customthm}

(Recall a field $K$ is \textit{bounded} if for every $n > 1$ there are finitely many Galois extensions $L/K$ with $[L : K] = n$.)

\bigskip
\subsection{Background} Recall specifically from \cite{haranjardenppc} the category of \textit{$G_{\Q_p}\!$-structures} and the following fact:

\begin{thm}
\emph{(\cite[Theorems 15.1 \& 15.3]{haranjardenppc}; cf.\ \emph{Theorem 15.4 ibid.})} If $K$ is a P$p$C field, then $\mathfrak{G}_K$ -- \emph{the absolute $G_{\Q_p}\!$-structure of $K$} -- is projective. Conversely, if $\mathfrak{G}$ is a projective $G_{\Q_p}\!$-structure, there exists a P$p$C field $K$ such that $\mathfrak{G} \iso \mathfrak{G}_K$. \bsn
\end{thm}

Again one can show ``projectivity'' -- we mean \cite[Definition 5.1]{haranjardenppc} -- coincides with ``projective in the category of $G_{\Q_p}\!$-structures''. See also the definition of a \emph{$G_{\Q_p}\!$-projective group} \cite[Definition 4.1]{haranjardenppc}, as these groups underlie projective $G_{\Q_p}\!$-structures exactly (\cite[Proposition 5.4]{haranjardenppc}). We have a corresponding \emph{Theorem \ref{genlv}/Corollary \ref{214}} as follows:

\begin{definition}
A field extension $L/K$ is \emph{totally $p$-adic} if the restriction map\footnote{See \cite[p.\ 180]{haranjardenppc} for the \emph{space of sites} $X(F) = X(F/F)$ of a field $F$.}\linebreak $\mbox{Res}: X(L) \rightarrow X(K)$ is surjective.
\end{definition}

This is to say we have `correctly' extended all the $p$-adic valuation data from $K$ to $L$. In the PRC setting we used \emph{totally real extensions} (on p.\ \pageref{totreal}) via {a similar characterisation}.

\begin{thm}\label{genlv3}
\emph{\cite[Theorem 15.3]{haranjardenppc}.} Let $\mathfrak{G}$ be a projective $G_{\Q_p}\!$-structure. Let $L/K$ be a Galois extension and $\pi : \mathfrak{G} \rightarrow \mathfrak{Gal}(L/K)$ be an epimorphism. Then there exists a \emph{totally $p$-adic} P$p$C extension $E/K$ such that $\mathfrak{G} \iso \mathfrak{G}_E$, and
\[
\begin{tikzcd}
\mathfrak{G} \arrow[r, "\iso"] \arrow[rd, "\pi", swap] & \mathfrak{G}_E \arrow[d, "\mbox{Res}_L"]\\
& \mathfrak{Gal}(L/K)
\end{tikzcd}
\]
is a commutative diagram. \bsn
\end{thm}

\noindent (Recall \cite[\S 10]{haranjardenppc} for the notation $\mathfrak{Gal}(L/K)$ here; cf.\ \cite[Proposition 10.7]{haranjardenppc}.) 

\begin{cor}\label{important}
Let $L/K$ be a Galois extension, $\mathfrak{G}$ a projective $G_{\Q_p}\!$-structure, and $\alpha : \mathfrak{G} \rightarrow \mathfrak{Gal}(L/K)$ an epimorphism. Then $K$ has a (totally $p$-adic) extension $E$ which is P$p$C, {linearly disjoint from $L$}, and there exists an isomorphism $\gamma : \mathfrak{G}_E \rightarrow \mathfrak{G}$ such that $\alpha \circ \gamma = Res_L$. \bsn
\end{cor}

As before, the Galois-theoretic information of a P$p$C field, along with some arithmetic information, completely determines that field's ($\mathcal{L}_r$-)theory:

\begin{thm}\label{thmforppcel}
Let $K_1, K_2$ be P$p$C fields, {separable} over a common subfield $E$. Then $K_1 \equiv_E K_2$ if and only if $K_1$ \&\ $K_2$ have {the same degree of imperfection,} there exists $\theta \in G_E$ such that $\theta(K_1 \cap E^s) = K_2 \cap E^s$, and if $S\Theta : S(G_{K_1 \cap E^s}) \rightarrow S(G_{K_2 \cap E^s})$ is the isomorphism induced by $\theta$, then the partial map $S\Theta : S(G_{K_1}) \rightarrow S(G_{K_2})$ with domain $S(G_{K_1 \cap E^s})$ is $\mathcal{L}_G$-elementary.
\end{thm}

\proofsketch
As in \emph{Theorem \ref{thmforprcel}}. Note this now requires a ``P$p$C embedding lemma'', i.e.\ \textit{if $K_1$ and $K_2$ are formally $p$-adic P$p$C fields that contain a common field $E$, and supposing that there exists a homeomorphism $\varphi: G_{K_1} \rightarrow G_{K_2}$ such that for every $\sigma \in G_{K_1}$, $\varphi(\sigma)|_{E^s} = \sigma|_{E^s}$; then $K_1 \equiv_E K_2$.} Such a lemma is provided by \cite[Lemma 5 \& Proposition 6]{kunzi}; cf.\ \cite[Proposition 10.4]{jardennew}.\bs

\begin{cor}\label{ppcax}
A P$p$C field $K$ with prime subfield $\mathbb{F}$ is axiomatised by the following first-order $(\mathcal{L}_r$-$)$axiom scheme:
\begin{enumerate}
    \item The characteristic and degree of imperfection of $K$;
    \item The P$p$C field axioms, \emph{\ttfamily P\emph{p}\hspace{-0.5mm}C};
    \item $\TTh^{alg}(K)$;
    \item $\TTh(S(G_K); \mathcal{L}_G(S(G_{K_0})))^*$, where $K_0 = K \cap \mathbb{F}^s$.\bs
\end{enumerate}
\end{cor}

Let us repeat \emph{Lemma \ref{lem2}} and prove there is a `counterfeit' coproduct in the category of $G_{\Q_p}\!$-structures, sufficiently similar to the real thing. Before this:

\begin{remark}\label{topologyremark}
{(cf.\ \cite[p.\ 150]{haranjardenppc}.)} Let $\mathfrak{G} = \langle G, d: X(\mathfrak{G}) \rightarrow \Hom(G_{\Q_p}, G) \rangle$ be a (weak) $G_{\Q_p}\!$-structure. Denote the set of closed subgroups of $G$ by $\Subg(G)$. Notice that $\Subg(G)$ can be viewed as a Boolean space: for open normal subgroups of $G$, equip $\Subg(G/N)$ with the discrete topology. By the compactness of $G$, $\Subg(G) \iso \varprojlim \Subg(G/N)$, hence $\Subg(G)$ is a profinite space under this topology. Furthermore, the map
$$\Im : \Hom(G_{\Q_p}, G) \rightarrow \Subg(G); \qquad \Im(\psi) = \psi(G_{\Q_p}),$$

is continuous, as it is the inverse limit of the continuous maps $\Im_N : \Hom(G_{\Q_p}, G/N) \rightarrow \Subg(G/N)$ of finite discrete spaces. This topology may be equivalently characterised as follows: called the \emph{strict topology}, it has basis
$$v(\Delta, N) = \{H \in \Subg(G) \mbox{ : } HN = \Delta N\}, \quad \mbox{where $\Delta \in \Subg(G)$ and $N \triangleleft G$ open.}$$

Recall from \cite{haranjardenppc} the notation $\mathcal{D}(G) = \{H \leq G \mbox{ : } H \iso G_{\Q_p}\}$ for a profinite group $G$. A necessary condition (\cite[Definition 4.1]{haranjardenppc}) for $G$ to be a $G_{\Q_p}$-projective group is for $\mathcal{D}(G)$ to be topologically closed in $\Subg(G)$, however this is \textit{always the case} {for $G_{\Q_p}\!$-structures, by a minor adaptation of} \cite[Lemma 3.5.1]{fehmthesis}:\\

\noindent \textbf{Claim.} \textit{Let $G$ be a profinite group. Then $\mathcal{D}(G)$ is topologically closed in $\Subg(G)$.}

\proof
Consider $H \in \Subg(G)$ such that $H \not\iso G_{\Q_p}$. As $G_{\Q_p}$ is finitely generated, by \cite[Lemma 1.3.2 (1)]{fehmthesis} there exists $H_0 \triangleleft H$ open such that $H/H_0$ is not a quotient of $G_{\Q_p}$. Fix $N \triangleleft G$ open with $N \cap H \leq H_0$ (which exists by \cite[Lemma 1.2.5 (a)]{friedjarden}). As $H/(N \cap H) \supseteq H/H_0$, $H/(N \cap H)$ is also not a quotient of $G_{\Q_p}\!$. If $H' \leq G$ and $H'N = HN$, then $H'/(N \cap H') \iso H'N/N \iso HN/N \iso H/(N\cap H)$. Hence $H'/(N \cap H')$ is also not a quotient of $G_{\Q_p}$, meaning $H' \not\iso G_{\Q_p}$. We have shown the ``strict'' neighbourhood $v(H, N) \subseteq \Subg(G)\setminus \mathcal{D}(G)$, hence $\mathcal{D}(G)$ is closed, as desired. \sq
\end{remark}

\begin{lem}\label{lem3}
Let $G$ be a $G_{\Q_p}$-projective group, and $E$ a projective profinite group. Then $G \star E$ is a $G_{\Q_p}\!$-projective group.
\end{lem}

\proof
Cf.\ \cite[\S 3.3]{fehmthesis}. Let the following be a finite $G_{\Q_p}$-embedding problem for $G \star E$:
\begin{equation}\label{firstdiag}
\begin{tikzcd}
& G\star E \arrow[d, "\varphi"]\\
B \arrow[r, "\alpha", twoheadrightarrow] & A
\end{tikzcd}    
\end{equation}
This diagram gives rise to the following two:
\[
\begin{tikzcd}
& G \arrow[d, "\varphi|_G"] &&&&E \arrow[d, "\varphi|_E"]\\
B \arrow[r, "\alpha", twoheadrightarrow] & A&&&B \arrow[r, "\alpha", twoheadrightarrow] & A
\end{tikzcd}
\]

The latter is solvable by $\gamma_E$, as $E$ is projective. We claim the former is in fact a $G_{\Q_p}$-embedding problem. Indeed, $\mathcal{D}(G) \hookrightarrow \mathcal{D}(G \star E)$ {as $G \leq G \star E$}, so if $H \in \mathcal{D}(G)$ there exists $\gamma_H : H \rightarrow B$ with $\alpha \circ \gamma_H = \varphi|_H = (\varphi|_G)|_H$. Therefore, as $G$ is $G_{\Q_p}$-projective, we conclude the former diagram is solvable by $\gamma_G$:
\[
\begin{tikzcd}
& G \arrow[d, "\varphi|_G"] \arrow[dl, dashrightarrow, swap, "\gamma_G"] &&&&E \arrow[d, "\varphi|_E"] \arrow[dl, dashrightarrow, swap, "\gamma_E"]\\
B \arrow[r, "\alpha", twoheadrightarrow] & A&&&B \arrow[r, "\alpha", twoheadrightarrow] & A
\end{tikzcd}
\]

By definition of the free product, there exists a morphism $\gamma : G \star E \rightarrow B$ {uniquely extending $\gamma_G$, $\gamma_E$, and} making (\ref{firstdiag}) commute. Finally, $\mathcal{D}(G\star E) = \{H \leq G\star E \mbox{ : } H \iso G_{\Q_p}\}$ is topologically closed in $\Subg(G\star E)$ by \emph{Remark \ref{topologyremark}}. This concludes the proof of the lemma. \bs

\subsection{Results}
Given a bounded P$p$C field $L$, we claim there does not exist an $\mathcal{L}_r$-sentence $\gamma$ such that for all P$p$C fields $F$ of the same characteristic and degree of imperfection as $L$, $F \models \gamma \iff F \equiv_{\mathcal{L}_r} L$. This is:

\begin{thm}\label{noppcfinax}
No bounded P$p$C field is finitely axiomatisable (even among the class of P$p$C fields).
\end{thm}

\proof
Assume for the purpose of contradiction that there is a formally $p$-adic finitely axiomatisable bounded P$p$C field $K$ (see \emph{Theorem \ref{pacfiniteaxiom}} for when $K$ is P$p$C but not formally $p$-adic). In order to characterise $K$ among P$p$C fields, by \emph{Corollary \ref{ppcax}} we need only finitely many axioms $\Sigma_1 \subset \Th^{alg}(K)$ and $\Sigma^*_2 \subset \Th(S(G_{K}); \mathcal{L}_G(S(G_{K_0})))^*$, where $\Sigma_2 \subset \Th(S(G_K); \mathcal{L}_G(S(G_{K_0})))$ is finite.

Let $\Lambda$ be the set of universal sentences of $\Sigma_1$, and let $\Kbad$ be the join of minimal Galois extensions $F/K$ within a fixed algebraic closure $\widetilde{K}$ of $K$, with $F \models \neg\lambda$ for $\lambda \in \Lambda$. Then $\Kbad/K$ is finite, $\Kbad \models \neg\Lambda$, and $\mbox{Res}: \mathfrak{G}_K \rightarrow \mathfrak{Gal}(\Kbad/K)$ is a cover (\cite[Proposition 10.7]{haranjardenppc}). As there is an epimorphism $G_K \twoheadrightarrow \Gal(\Kbad/K)$ -- apply the forgetful functor to Res -- there is an embedding of $\mathcal{L}_G$-structures $S(\Gal(\Kbad/K)) \hookrightarrow S(G_K)$.

Let $\overline{a} \in S(G_{K_0}) \subset S(G_K)$ be a finite tuple of elements such that $\Sigma_2$ is a set of finitely many $\mathcal{L}_G(\overline{a})$-sentences. Fix $n_k \in \N$ such that $S_1, \dots, S_{n_k}$ is the smallest consecutive sequence of sorts involving the sentences of $\Sigma_2$. Let $\widehat{p}$ be the smallest prime larger than $n_k + |\Gal(\Kbad/K)|$ and  $\mathbb{F}_1(q)$ be the free pro-$q$ group on one generator, where $q$ is a prime larger than $\widetilde{p}$. This is a projective profinite group. Consider $G_K \star \mathbb{F}_1(q)$; by \emph{Lemma \ref{lem3}} this is a $G_{\Q_p}\!$-projective group, and $S(G_K \star \mathbb{F}_1(q)) \models \Sigma_2$ by \emph{Lemma \ref{addedsea}} (with $\widetilde{G_{\Gamma}}$ replaced by $\mathbb{F}_1(q)$). However the $\mathcal{L}_G$-theories of $S(G_K \star \mathbb{F}_1(q))$ and $S(G_K)$ differ: indeed, as $K$ is bounded, $G_K$ is small (as a profinite group) hence there are $1 \leq n_q < \infty$ many open normal subgroups of $G_K$ of index $\leq q$. Consider the sentence
$$\varphi_{n_q}: \qquad \exists N_0, \dots, N_{n_q} \in S_q \, \left(\bigwedge_{\substack{i \neq j\\0 \leq i,j \leq n_q}} \neg (N_i \leq N_j \land N_j \leq N_i) \right),$$

where $S_q$ denotes sort $q$, $N_i = g_iM_i$ is a coset of an open $M_i \triangleleft G_K$ of index $\leq q$, and recall $N_i \leq N_j \iff M_i \subseteq M_j$. This sentence expresses ``there are $n_q+1$ distinct open normal subgroups of index $\leq q$'', and clearly $S(G_K) \not\models \varphi_{n_q}$ while $S(G_K \star \mathbb{F}_1(q)) \models \varphi_{n_q}$.

Let $\mathfrak{G}_K \Diamond \mathbb{F}_1(q)$ be a \textit{projective $G_{\Q_p}\!$-structure} with underlying group $G_K \star \widetilde{G_{\Gamma}}$ (this is \cite[Proposition 5.4]{haranjardenppc}). By \emph{Theorem \ref{genlv3}} (setting $L = K$) there exist totally $p$-adic P$p$C fields $E_1, E_2 \supseteq K$ such that $\mathfrak{G}_{E_1} \iso \mathfrak{G}_K$ and $\mathfrak{G}_{E_2} \iso \mathfrak{G}_K \Diamond \mathbb{F}_1(q)$. As there are epimorphisms $G_{E_1} \iso G_K \twoheadrightarrow \mbox{Gal}(\Kbad/K)$, $G_{E_2} \iso G_K \star \widetilde{G_{\Gamma}} \twoheadrightarrow G_K \twoheadrightarrow \mbox{Gal}(\Kbad/K)$, $E_1 \cap \Kbad = E_2 \cap \Kbad = K$, therefore $E_1, E_2 \models \Sigma_1$. We conclude $E_1, E_2 \models \Sigma_1 \cup \Sigma_2^*$ however $E_1 \not\equiv E_2$ as $S(G_{E_1}) \not\equiv_{\mathcal{L}_G} S(G_{E_2})$; a contradiction, as required. \bs

\begin{remark}
Again, there is a classification-theoretic connection: every bounded P$p$C field has an NTP${}_2$ $\mathcal{L}_r$-theory \cite[Corollary 8.6]{mont}, and it is conjectured (\cite[Conjecture 5.1]{cks}) that NTP${}_2$ P$p$C fields are bounded. \sq
\end{remark}

\begin{remark}\label{lackofppc}
We are unable to prove formally $p$-adic P$p$C fields are finitely undecidable -- indeed, difficulties arise adapting the previous proofs to the $p$-adic context, simply because $G_{\Q_p}\!$ has a more complicated $\mathcal{L}_G$-structure theory than $\Z\!/2\Z$ ($= G_{\R}$) or $\{1\}$ ($=G_{\C}$). 

The first fundamental issue is that, if $\mathcal{C}$ is a full formation of finite $P$-groups where $P$ is a finite set of consecutive primes $\{2, 3, \dots, \widehat{p}\}$, it is not true that the maximal pro-$\mathcal{C}$ quotient of a $G_{\Q_p}\!$-projective group remains $G_{\Q_p}\!$-projective. Indeed, the absolute Galois group of a formally $p$-adic P$p$C field $K$ admits embeddings $G_{\Q_p} \hookrightarrow G_K$ by \cite[Lemma 5.3]{haranjardenppc}; impossible if $G_K$ is pro-$\mathcal{C}$.

The second fundamental issue is the following: let $K$ be a PAC (resp.\ PRC) field, $\Sigma \subset \Th(K; \mathcal{L}_r)$ a finite subtheory, and $\Gamma$ be a nontrivial graph. In \emph{Theorem \ref{thisolthm}} (resp.\ \emph{Theorem \ref{mainthm2}}), the field $K_{\Gamma}$ constructed has absolute Galois group $H \star \widetilde{G_{\Gamma}}$ where $H$ is projective (resp.\ real projective) and has no open normal subgroups of index $n_{\Sigma} = |D_{\Sigma}|$ or $|W_{\Sigma}|$. If $H$ is the absolute Galois group of a formally $p$-adic P$p$C field, we cannot be guaranteed that any index $n_{\Sigma}$ open normal subgroup of $H \star \widetilde{G_{\Gamma}}$ contains $H$. Consequently, although $\Gamma \hookrightarrow \Gamma_{H\star\widetilde{G_{\Gamma}}}$ ($= \Gamma_{G_{K_{\Gamma}}}$, interpretable in $K_{\Gamma}$ by \cite[Lemma 4.1]{efrat}) we cannot guarantee these graphs are isomorphic. \sq
\end{remark}

\section*{Acknowledgements}
The author extends his thanks to Professor Ehud Hrushovski, Professor Jochen \linebreak Koenigsmann, and Professor Arno Fehm for their time, comments, and efforts. I am grateful to Mikhail Blinov for his translation of \cite[Chapter 5, \S 1.4]{ershovrussian} for me (any errors in presentation are, of course, my own).

\bigskip
\bibliography{refs}   
\bibliographystyle{acm} 
\end{document}